\documentclass[11pt]{amsart}
\usepackage{amsopn}
\usepackage{amsmath,amsthm,amssymb}

\newcommand{\nc}{\newcommand}

\nc{\n}{\noindent}
\nc{\Id}{\mathbf{1}}
\nc{\vs}{\vspace{8pt}}
\nc{\alt}{\raise1pt\hbox{$\bigwedge$}}
\nc{\ncp}{\nabla^\mathrm{CP}}
\nc{\nhc}{\nabla^\mathrm{HC}}
\nc{\wncp}{\widehat\nabla^\mathrm{CP}}
\nc{\ct}{\cos_\theta}
\nc{\st}{\sin_\theta}
\nc{\ctt}{\cos_{\theta/2}}
\nc{\stt}{\sin_{\theta/2}}
\nc{\ft}{\hbox{$\frac12$}}
\nc{\Aff}{\mathit{Aff}}

\nc{\noi}{\noindent}
\nc{\vg}{\mathfrak{v} }
\nc{\wg}{\mathfrak{w} }
\nc{\zg}{\mathfrak{z} }
\nc{\ngo}{\mathfrak{n} }
\nc{\kg}{\mathfrak{k} }
\nc{\mg}{\mathfrak{m} }
\nc{\bg}{\mathfrak{b} }
\nc{\ggo}{\mathfrak{g} }
\nc{\ggob}{\overline{\mathfrak{g}} }
\nc{\sog}{\mathfrak{so} }
\nc{\sug}{\mathfrak{su} }
\nc{\spg}{\mathfrak{sp} }
\nc{\slg}{\mathfrak{sl} }
\nc{\glg}{\mathfrak{gl} }
\nc{\cg}{\mathfrak{c} }
\nc{\hg}{\mathfrak{h} }
\nc{\tg}{\mathfrak{t} }
\nc{\ug}{\mathfrak{u} }
\nc{\dg}{\mathfrak{d} }
\nc{\ag}{\mathfrak{a} }
\nc{\pg}{\mathfrak{p} }
\nc{\sg}{\mathfrak{s} }
\nc{\rg}{\mathfrak{r} }
\nc{\aff}{\mathfrak{aff}}
\nc{\pca}{\mathcal{P}}
\nc{\nca}{\mathcal{N}}

\nc{\vp}{\varphi}
\nc{\ddt}{\frac{{\rm d}}{{\rm d}t}}

\nc{\SO}{{\tt SO}}
\nc{\Spe}{{\tt Sp}}
\nc{\SL}{{\tt SL}}
\nc{\SU}{{\tt SU}}
\nc{\Or}{{\tt O}}
\nc{\U}{{\tt U}}
\nc{\GL}{{\tt GL}}
\nc{\Se}{{\tt S}}
\nc{\CL}{{\tt CL}}
\nc{\Spin}{{\tt Spin}}
\nc{\Pin}{{\tt Pin}}

\nc{\RR}{{\mathbb R}}
\nc{\HH}{{\mathbb H}}
\nc{\CC}{{\mathbb C}}
\nc{\ZZ}{{\mathbb Z}}
\nc{\FF}{{\mathbb F}}
\nc{\NN}{{\mathbb N}}
\nc{\GG}{{\mathbb G}}
\nc{\JJ}{{\mathbb J}}
\nc{\II}{{\mathbb I}}
\nc{\KK}{{\mathbb K}}
\nc{\DD}{{\mathbb D}}

\nc{\ad}{\operatorname{ad}}
\nc{\Ad}{\operatorname{Ad}}
\nc{\rank}{\operatorname{rank}}
\nc{\Irr}{\operatorname{Irr}}
\nc{\End}{\operatorname{End}}
\nc{\Aut}{\operatorname{Aut}}
\nc{\Inn}{\operatorname{Inn}}
\nc{\Der}{\operatorname{Der}}
\nc{\Ker}{\operatorname{Ker}}
\nc{\Iso}{\operatorname{I}}
\nc{\Le}{\operatorname{L}}
\nc{\tr}{\operatorname{tr}}
\nc{\dif}{\operatorname{d}\!}
\nc{\sen}{\operatorname{sen}}
\nc{\modu}{\operatorname{mod}}
\nc{\Ric}{\operatorname{R}}
\nc{\Sym}{\operatorname{Sym}}
\nc{\sca}{\operatorname{sc}}
\nc{\scalar}{{\sf s}}
\nc{\grad}{\operatorname{grad}}
\nc{\ricci}{\operatorname{r}}
\nc{\riccin}{\operatorname{Ric}}
\nc{\Lie}{\operatorname{L}}
\nc{\tang}{\operatorname{T}}
\nc{\tm}{\operatorname{TM}}

\theoremstyle{plain}
\newtheorem{thm}{Theorem}[section]
\newtheorem{prop}[thm]{Proposition}

\newtheorem{lem}[thm]{Lemma}

\theoremstyle{definition}
\newtheorem{defi}[thm]{Definition}

\theoremstyle{remark}
\newtheorem*{rem}{Remark}
\newtheorem*{rems}{Remarks}

\newcommand{\ri}{{\rm (i)}}
\newcommand{\rii}{{\rm (ii)}}
\newcommand{\riii}{{\rm (iii)}}
\newcommand{\riv}{{\rm (iv)}}

\title{Hypersymplectic four-dimensional Lie algebras}

\author[A.~Andrada]{Adri\'an Andrada}

\address{CIEM, FaMAF, Universidad Nacional de C\'ordoba, Ciudad Universitaria,
(5000) C\'ordoba, Argentina}
\email{andrada@mate.uncor.edu}

\begin{document}

\maketitle

\n\textbf{Abstract.} A study is made of real Lie algebras admitting a hypersymplectic structure, and we provide a method to construct such hypersymplectic Lie algebras. We use this method in order to obtain the classification of all hypersymplectic structures on four-dimensional Lie algebras, and we describe the associated metrics on the corresponding Lie groups.

\smallbreak\n\textbf{MSC.} 17B60, 53C15, 53C30, 53C50. 

\section{Introduction}

A hypersymplectic structure on a manifold is a complex product
structure, i.e. a pair $\{J,E\}$ of a complex structure and a
product structure that anticommute, together with a compatible
metric such that the associated 2-forms are closed. This notion is
similar to that of a hyperk\"ahler structure, where the base
manifold carries a hypercomplex structure, i.e. a pair
$\{J_1,J_2\}$ of anticommuting complex structures.

Hypersymplectic structures were introduced by N.~Hitchin in
\cite{H1}, and they are also referred to as {\em neutral
hyperk\"ahler} structures in \cite{Kam} and as {\em
parahyperk\"ahler} structures in \cite{V}. Hypersymplectic
structures on manifolds have become an important subject of study
lately, due mainly to its applications in theoretical physics
(specially in dimension 4). See for instance \cite{BGPPR}, where
there is a discussion on the relationship between hypersymplectic
metrics and the $N=2$ superstring. Hypersymplectic metrics on a
manifold are Ricci-flat and the associated holonomy group is
contained in the real symplectic group.

In \cite{Kam}, H.~Kamada determines the compact complex surfaces
which admit hypersymplectic structures. These complex surfaces are
either complex tori or primary Kodaira surfaces; Kamada also shows
when the hypersymplectic metrics on these surfaces are flat. In
\cite{FPPS}, examples of (non flat) hypersymplectic structures are
given on Kodaira manifolds, which are special compact quotients of
2-step nilpotent Lie groups. These hypersymplectic structures are
not invariant by the nilpotent Lie group.

The main goal of this paper is to give the classification, up to
equivalence, of all left-invariant hypersymplectic structures on
4-dimensional Lie groups. These Lie groups will provide examples
of hypersymplectic structures in non compact manifolds, since
their underlying differentiable manifolds are diffeomorphic to
$\RR^4$. In order to perform this classification, we begin in \S 3
the study of hypersymplectic structures on real Lie algebras. We
obtain that, associated to a hypersymplectic structure $\{J,E,g\}$
on a Lie algebra $\ggo$, there are two 3-tuples
$(\ggo_+,\nabla^+,\omega_+),\,(\ggo_-,\nabla^-,\omega_-)$, where
$\ggo_{\pm}$ are Lie subalgebras of $\ggo$ such that
$\ggo=\ggo_+\oplus\ggo_-$ and $\ggo_-=J\ggo_+$, $\nabla^{\pm}$ is
a flat torsion-free connection on $\ggo_{\pm}$ and $\omega_{\pm}$
is a symplectic form on $\ggo_{\pm}$ such that
$\omega_+(x,y)=\omega_-(Jx,Jy)$ for all $x,y\in\ggo_+$ and
$\nabla^{\pm}\omega_{\pm}=0$. Conversely, we show that, in certain
cases, given two 3-tuples $(\ug,\nabla,\omega)$ and
$(\vg,\nabla',\omega')$ satisfying the same conditions as above,
we can obtain a hypersymplectic structure on $\ug\oplus\vg$
(direct sum of vector spaces). This result will be used in the
4-dimensional case. We also deal with equivalences between
hypersymplectic structures.

Next, in \S 4, we give the first steps in order to achieve the
classification mentioned above, namely, we determine the flat
torsion-free connections on the 2-dimensional Lie algebras which
are compatible with a symplectic form and obtain their equivalence
classes.

In \S 5, we prove our main result which states that, aside from
the abelian Lie algebra, there are only three Lie algebras which
admit a hypersymplectic structure. One of them is a central
extension of the $3$-dimensional Heisenberg algebra $\hg_3$; the
second one is an extension of $\RR^3$ and the third one is an
extension of $\hg_3$. We also parameterize the underlying complex
product structures. In \S 6, we show that a complex product
structure on a Lie algebra admits at most one compatible metric
(up to a multiplicative constant), and we determine the
hypersymplectic metrics for each of the complex product structures
in the Lie algebras obtained previously. We point out the cases
when these metrics are flat and/or complete. As an illustration we
exhibit the following examples of hypersymplectic metrics on
$\RR^4$ with canonical global coordinates $t,x,y,z$:
\begin{enumerate}
\item[\ri] $g=\dif t^2+\dif x^2-\dif y^2-\dif z^2$ (flat and
complete). \item[\rii] $g=e^{-t}\dif t(\dif z-\ft x\dif y+\ft
y\dif x)+e^{-t}\dif x\dif y$ (flat but not complete).
\item[\riii] $g=e^t\dif t\dif z+ e^{2t}\dif z^2-\dif x\dif y+
e^{2t}\dif y^2$ (neither flat nor complete).
\end{enumerate}
The three metrics above are hence not isometric.

\vs

\noindent{\textbf{Acknowledgements.} This work is part of the author's doctoral thesis under the guidance of Prof. Isabel G. Dotti., to whom the author would like to thank for supervision and invaluable comments.  The author was supported by grants from CONICET, FONCYT and SECYT-UNC (Argentina).}

\

\section{Preliminaries}

We start recalling some definitions which will be used throughout
this work. All Lie algebras will be finite dimensional and defined
over $\RR$, unless explicitly stated.

For an arbitrary connection $\nabla$ on (the tangent bundle of) a
manifold $M$, the torsion and curvature tensor fields $T$ and $R$
are defined by \begin{gather*}
T(X,Y)=\nabla_XY-\nabla_YX-[X,Y]\\R(X,Y)=
\nabla_X\nabla_Y-\nabla_Y\nabla_X-\nabla_{[X,Y]} \end{gather*} for
$X,Y$ smooth vector fields on $M$. The connection is called {\em
torsion-free} when $T=0$, and {\em flat} when $R=0$.

Let $G$ be a Lie group with Lie algebra $\ggo$ and suppose that
$G$ admits a left-invariant connection $\nabla$. This means that
if $X,Y\in\ggo$ are two left-invariant vector fields on $G$ then
$\nabla_XY\in\ggo$ is also left-invariant. Accordingly, one may
define a connection on a Lie algebra $\ggo$ to be merely a
$\ggo$-valued bilinear form $\ggo\times\ggo \longrightarrow\ggo$.
One can speak of the torsion and curvature of such a connection
using the formulae above, with brackets determined by the
structure of $\ggo$. It is known that the completeness of the
left-invariant connection $\nabla$ on $G$ can be studied by
considering simply the corresponding connection on the Lie algebra
$\ggo$. Indeed, the connection $\nabla$ on $G$ will be
(geodesically) complete if and only if the differential equation
on $\ggo$
\begin{equation}\label{complete}
\dot{x}(t)=-\nabla_{x(t)}x(t)
\end{equation}
admits solutions $x(t)\in\ggo$ defined for all $t\in\RR$ (see for instance
\cite{G}, where only metric connections are
considered).

\vs

We also recall the definition of complex structures and product
structures on a Lie algebra, which are also modelled on the
corresponding notions for smooth manifolds.

An {\em almost complex structure} on a Lie algebra $\ggo$ is a linear
endomorphism $J:\ggo \longrightarrow \ggo$ satisfying $J^2=-\Id$. If $J$
satisfies the condition
\begin{equation} \label{integrable}
J[X,Y]=[JX,Y]+[X,JY]+J[JX,JY] \quad \text{for all } X,Y\in \ggo,
\end{equation}
we will say that $J$ is {\em integrable} and we will call it a
{\em complex structure} on $\ggo$. Note that the dimension of a
Lie algebra carrying an almost complex structure must be even. We
recall that a {\em hypercomplex structure} on the Lie algebra
$\ggo$ is a pair $\{J_1,J_2\}$ of complex structures on $\ggo$
such that $J_1J_2=-J_2J_1$. The dimension of a hypercomplex Lie
algebra is a multiple of 4.

Next, an {\em almost product structure} on $\ggo$ is a linear
endomorphism $E:\ggo \longrightarrow \ggo$ satisfying $E^2=\Id$
(and not equal to $\pm\Id$). It is said to be {\em integrable} if
\begin{equation} \label{integrable2} E[X,Y]=[EX,Y]+[X,EY]-E[EX,EY]
\quad \text{for all } X,Y\in \ggo.  \end{equation} An integrable
almost product structure will be called a {\em product structure}.
If $\ggo_\pm$ is the eigenspace of $\ggo$ associated to the
eigenvalue $\pm1$ of $E$, then the integrability of $E$ is
equivalent to the fact of $\ggo_{\pm}$ being Lie subalgebras of
$\ggo$. If $\dim \ggo_+=\dim\ggo_-$, the product structure $E$ is
called a {\em paracomplex structure} \cite{KK,L}. In this case,
$\ggo$ also has even dimension.

\

An appropriate combination of these two structures on Lie algebras
is called a complex product structure, and its definition is given
below. This new structure is similar to a hypercomplex structure,
where one of the complex structures has been replaced by a product
structure.

\smallskip

\begin{defi} 
A {\em complex product structure} on the Lie algebra $\ggo$ is a pair
$\{J,\,E\}$ of a complex structure $J$ and a product structure $E$ satisfying
$JE=-EJ$.
\end{defi}

Complex product structures on Lie algebras have been studied in
\cite{AS}, from where we recall some of their main properties. The
condition $JE=-EJ$ implies that $J$ is an isomorphism (as vector
spaces) between the eigenspaces $\ggo_+$ and $\ggo_-$
corresponding to the eigenvalues $+1$ and $-1$ of $E$,
respectively; thus, $E$ is in fact a paracomplex structure on
$\ggo$. Every complex product structure on $\ggo$ has therefore an
associated {\em double Lie algebra} $(\ggo,\ggo_+,\ggo_-)$, i.e.,
$\ggo_+$ and $\ggo_-$ are Lie subalgebras of $\ggo$ such that
$\ggo=\ggo_+\oplus\ggo_-$ (direct sum of vector spaces) and
$\ggo_-=J\ggo_+$, where $E|_{\ggo_+}=\Id,\,E|_{\ggo_+}=-\Id$. We
note that the dimension of a Lie algebra with a complex product
structure is even but it need not be a multiple of $4$.

\medskip

The complex product structure $\{J,E\}$ on $\ggo$ determines
uniquely a torsion-free connection $\ncp$ on $\ggo$ such that
$\ncp J=\ncp E=0$, where this equations mean that
\[ \ncp_x Jy=J\ncp _xy,\qquad \ncp_x Ey=E\ncp _xy\]
for all $x,y\in\ggo$. As a consequence, we note that
$\ncp_xy\in\ggo_{\pm}$ for any $x\in\ggo$ and $y\in\ggo_{\pm}$.
Take now $x\in\ggo_+,y\in\ggo_-$. Since $\ncp$ has no torsion, we
obtain that \begin{equation}\label{split}
[x,y]=-\ncp_yx+\ncp_xy\in\ggo_+\oplus\ggo_-\end{equation} is the
decomposition of $[x,y]$ into components, according to the
splitting $\ggo=\ggo_+\oplus\ggo_-$. The connection $\ncp$
restricts to {\em flat} torsion-free connections on $\ggo_+$ and
$\ggo_-$, say $\nabla^+$ and $\nabla^-$, respectively. Thus, we
have that $\ncp$ is flat if and only if $R(x_+,x_-)=0$ for all
$x_+\in\ggo_+,x_-\in\ggo_-$, where $R$ is the curvature of $\ncp$.
We recall that flat torsion-free connections on a Lie algebra are
also known as ``left-symmetric algebra" (LSA) structures.

\section{Hypersymplectic structures on Lie algebras}

We have observed the close resemblance between complex product
structures and hypercomplex structures on Lie algebras. Bearing
this similarity in mind, we study in this section a special kind
of metrics on a Lie algebra with a complex product structure, just
as hyperk\"ahler metrics appear in the context of hypercomplex
structures.

Let $\ggo$ be a Lie algebra endowed with a complex product
structure $\{J,E\}$ and let $g$ be a metric on $\ggo$, i.e., $g$
is a non degenerate symmetric bilinear form
$g:\ggo\times\ggo\longrightarrow\RR$. We will say that $g$ is {\em
compatible} with the complex product structure if
\begin{equation}\label{metric}
g(Jx,Jy)=g(x,y),\quad g(Ex,Ey)=-g(x,y)
\end{equation}
for all $x,y\in\ggo$.

Let $(\ggo,\ggo_+,\ggo_-)$ denote the double Lie algebra
associated to the complex product structure $\{J,E\}$, where
$\ggo_-=J\ggo_+$ and let $g$ be a compatible metric. Then the
subalgebras $\ggo_+$ and $\ggo_-$ are isotropic subspaces of
$\ggo$ with respect to $g$, i.e.,
\[ g(\ggo_+,\ggo_+)=0,\qquad g(\ggo_-,\ggo_-)=0.\]
For if $x,y\in\ggo_+$, we have that $g(Ex,Ey)=g(x,y)$ because of
the definition of $\ggo_+$. But $g(Ex,Ey)=-g(x,y)$ due to
(\ref{metric}). Hence $g(x,y)=0$ for $x,y\in\ggo_+$ and equally
for $x,y\in\ggo_-$. From this it is clear that
$\ggo_+^{\perp}=\ggo_+$ and $\ggo_-^{\perp}=\ggo_-$ and also that
the signature of $g$ is $(m,m)$, where $\dim\ggo=2m$.

Let us now define the following bilinear forms on $\ggo$:
\begin{equation}\label{forms} \omega_1(x,y)=g(Jx,y),\quad \omega_2(x,y)=g(Ex,y),\quad
\omega_3(x,y)=g(JEx,y)\end{equation} for $x,y\in\ggo$. Using
(\ref{metric}), it is readily verified that these forms are in
fact skew-symmetric, so that $\omega_i\in\alt^2\ggo^*$ for
$i=1,2,3$. Note that these forms are non degenerate, since $g$ is
non degenerate and $J$ and $E$ are isomorphisms. In the following
result we show the existing relationships between these 2-forms on
$\ggo$ and the decomposition of this Lie algebra induced by the
product structure $E$.

\begin{lem}\label{lema}
The 2-forms $\omega_i,\,i=1,2,3,$ on $\ggo$ satisfy the following
properties:
\begin{enumerate}
\item[\ri] $\omega_1(x,y)=\omega_1(Jx,Jy)=\omega_1(Ex,Ey)$ for all
$x,y\in\ggo$, whence $\omega_1(x,y)=0$ for
$x\in\ggo_+,\,y\in\ggo_-$. \item[\rii]
$\omega_2(x,y)=-\omega_2(Jx,Jy)=-\omega_2(Ex,Ey)$ for any
$x,y\in\ggo$, whence $\omega_2(x,y)=0$ for $x,y\in\ggo_+$ or
$x,y\in\ggo_-$. \item[\riii]
$\omega_3(x,y)=-\omega_3(Jx,Jy)=\omega_3(Ex,Ey)$ for all
$x,y\in\ggo$, whence $\omega_3(x,y)=0$ for
$x\in\ggo_+,\,y\in\ggo_-$.
\end{enumerate}
\end{lem}

\begin{proof}
The proof is straightforward.
\end{proof}

Let $\omega_+$ and $\omega_-$ denote the restriction of $\omega_1$
to $\ggo_+$ and $\ggo_-$, respectively. From $\ri$ of the previous
lemma and the fact that $\omega_1$ is non degenerate, it is easy
to see that both $\omega_+$ and $\omega_-$ are non degenerate.
Hence, $m=\dim\ggo_+=\dim\ggo_-$ must be an even number, say
$m=2n$, and therefore $\dim\ggo=4n$ and the signature of $g$ is
$(2n,2n)$.

From Lemma \ref{lema} $\ri$, we obtain that
\begin{equation}\label{equiv-simp} \omega_+(x,y)=\omega_-(Jx,Jy) \end{equation}
for all $x,y\in\ggo_+$. We shall show next that the forms
$\omega_1,\,\omega_2$ and $\omega_3$ can be written in terms
exclusively of $\omega_+$. In fact, we will show that
\begin{align}
\omega_1(x+Jx',y+Jy') & =\omega_+(x,y)+\omega_+(x',y'),\label{om1}\\
\omega_2(x+Jx',y+Jy') & =-\omega_+(x,y')+\omega_+(y,x'),\label{om2}\\
\omega_3(x+Jx',y+Jy') & =\omega_+(x,y)-\omega_+(x',y')\label{om3}
\end{align}
for all $x,y,x',y'\in \ggo_+$. Indeed, to prove (\ref{om1}) we
compute
\begin{align*}
\omega_1(x+Jx',y+Jy') & =\omega_1(x,y)+\omega_1(Jx',Jy')\\
& =\omega_1(x,y)+\omega_1(x',y')\\
& =\omega_+(x,y)+\omega_+(x',y')
\end{align*}
using Lemma \ref{lema} $\ri$. Next, we note that, for any
$u,v\in\ggo_+$ we have
\[ \omega_2(u,Jv)=g(Eu,Jv)=-g(JEu,v)=-g(Ju,v)=-\omega_1(u,v).\]
Now, in order to prove (\ref{om2}) we compute \begin{align*}
\omega_2(x+Jx',y+Jy') & = \omega_2(x,Jy')+\omega_2(Jx',y)\\
& = -\omega_1(x,y')+\omega_1(y,x')\\
& = -\omega_+(x,y')+\omega_+(y,x') \end{align*} because of Lemma
\ref{lema} $\rii$.  Finally, to verify that (\ref{om3}) holds, we
note first that from (\ref{forms}) we see that
$\omega_3(u,v)=\omega_1(u,v)$ for $u,v\in\ggo_+$; as a
consequence, we have
\begin{align*}
\omega_3(x+Jx',y+Jy') & = \omega_3(x,y)+\omega_3(Jx',Jy') \\
& = \omega_3(x,y)-\omega_3(x',y')\\
& = \omega_1(x,y)-\omega_1(x',y')\\
& = \omega_+(x,y)-\omega_+(x',y'), \end{align*} using Lemma
\ref{lema} $\riii$.

\

Let us recall now that given a 2-form $\omega$ on a Lie algebra
$\ggo$, there is an associated 3-form $\dif\omega\in\alt^3\ggo^*$
given by
\[ (\dif\omega)(x,y,z)=-\omega([x,y],z)+\omega([x,z],y)-\omega([y,z],x) \]
for all $x,y,z\in\ggo$. The 2-form $\omega$ is called {\em closed}
if $\dif\omega=0$; if $\omega$ is non degenerate and closed, it is
called a {\em symplectic form} on $\ggo$.

Naturally, we are mainly interested in the case when all of the
2-forms given in (\ref{forms}) are closed and hence symplectic.
We introduce therefore the following definition, equivalent to the
one given by N.~Hitchin in \cite{H1}.

\

\begin{defi}
Let $\{J,E\}$ be a complex product structure on the Lie algebra
$\ggo$ and let $g$ be a metric on $\ggo$ compatible with
$\{J,E\}$. If the 2-forms $\omega_i\in\alt^2\ggo^*$ defined in
(\ref{forms}) are closed, we will say that $\{J,E,g\}$ is a {\em
hypersymplectic structure} on $\ggo$. The Lie algebra $\ggo$ will
be referred to as a {\em hypersymplectic Lie algebra} and $g$ will
be called a {\em hypersymplectic metric}.
\end{defi}

\

The surprising fact is that if one of the 2-forms $\omega_1$ or
$\omega_3$ is closed, then all three of these 2-forms are closed,
as the following result shows.

\begin{prop}\label{omegas}
Let $\{J,E\}$ be a complex product structure on $\ggo$ with
associated double Lie algebra $(\ggo,\ggo_+,\ggo_-)$. Let
$\nabla^+$ and $\nabla^-$ denote the flat torsion-free connections
on $\ggo_+$ and $\ggo_-$ induced by $\ncp$. Suppose $g$ is a
compatible metric on $\ggo$ and let $\omega_i,\,i=1,2,3$, be the
2-forms on $\ggo$ given by (\ref{forms}) and $\omega_+$ and
$\omega_-$ be as above. Then the following statements are
equivalent:
\begin{enumerate}
\item[\ri] $\omega_1$ is closed; \item[\rii] $\omega_3$ is closed;
\item[\riii] $\nabla^+\omega_+=0$ and $\nabla^-\omega_-=0$.
\end{enumerate}
Furthermore, if one of the conditions above holds, then
\begin{enumerate}
\item[\riv] $\omega_2$ is closed.
\end{enumerate}
\end{prop}

\begin{rem} We recall that a 2-form $\omega$ on a Lie algebra $\hg$ is
parallel with respect to a connection $\nabla$ on $\hg$, i.e.
$\nabla\omega=0$,  if the condition
\[ \omega(\nabla_xy,z)=\omega(\nabla_xz,y)\]
holds for all $x,y,z\in\hg$.\end{rem}

\begin{proof} $\ri\Leftrightarrow\riii$
Let us suppose first that $\ri$ holds. For $x,y\in\ggo_+$ and
$z=Ju$ with $u\in\ggo_+$ we have that
\begin{align}
(\dif\omega_1)(x,y,Ju) & =\omega_1([x,Ju],y)-\omega_1([y,Ju],x)\label{nose}\\
& =-\omega_1(\ncp_{Ju}x,y)+\omega_1(\ncp_{Ju}y,x)\nonumber\\
& =-\omega_1(\ncp_{Ju}Jx,Jy)+\omega_1(\ncp_{Ju}Jy,Jx)\nonumber\\
&
=-\omega_-(\nabla^-_{Ju}Jx,Jy)+\omega_-(\nabla^-_{Ju}Jy,Jx),\nonumber
\end{align} using (\ref{split}). As $\dif\omega_1=0$, we obtain that
$\nabla^-\omega_-=0$. If we consider now $x\in\ggo_+$ and
$y=Jv,\,z=Ju$ with $u,v\in\ggo_+$, we have that
\begin{align}
(\dif\omega_1)(x,Jv,Ju) &=-\omega_1([x,Jv],Ju)+\omega_1([x,Ju],Jv)\label{nose2}\\
& =-\omega_1(\ncp_xJv,Ju)+\omega_1(\ncp_xJu,Jv)\nonumber\\
& =-\omega_1(\ncp_xv,u)+\omega_1(\ncp_xu,v)\nonumber\\
& =-\omega_+(\nabla^+_xv,u)+\omega_+(\nabla^+_xu,v),\nonumber
\end{align} using again (\ref{split}). As $\dif\omega_1=0$,
we obtain that $\nabla^+\omega_+=0$. Thus, $\riii$ holds.

Conversely, let us suppose that $\riii$ holds. We note first that,
as $\nabla^+$ and $\nabla^-$ are torsion-free, one obtains that
$\dif\omega_+=0$ and $\dif\omega_-=0$. Suppose that
$x,y,z\in\ggo_+$. Then $\dif\omega_1(x,y,z)=\dif\omega_+(x,y,z)=0$
since $\omega_+$ is closed. Similarly, for $x,y,z\in\ggo_-$, we
have $\dif\omega_1(x,y,z)=\dif\omega_-(x,y,z)=0$ since $\omega_-$
is closed. Now, if $x,y\in\ggo_+$ and $z=Ju$ with $u\in\ggo_+$,
from equations (\ref{nose}) and $\nabla^-\omega_-=0$, we have that
$(\dif\omega_1)(x,y,Ju)=0$. Next, if $x\in\ggo_+$ and
$y=Jv,\,z=Ju$ with $u,v\in\ggo_+$, from equations (\ref{nose2})
and $\nabla^+\omega_+=0$, we have that
$(\dif\omega_1)(x,Jv,Ju)=0$. Therefore, $\dif\omega_1=0$.

\smallskip

$\rii\Leftrightarrow\riii$ The proof is similar to the proof of $\ri\Leftrightarrow\riii$.

\smallskip

$\riii\Rightarrow\riv$ If $x,y,z\in\ggo_+$ or $x,y,z\in\ggo_-$,
then $(\dif\omega_2)(x,y,z)=0$, because of Lemma \ref{lema}
$\rii$. If $x,y\in\ggo_+$ and $z=Ju$ with $u\in\ggo_+$, then
\begin{align*} (\dif\omega_2)(x,y,Ju) & =
-\omega_2([x,y],Ju)+\omega_2([x,Ju],y)-\omega_2([y,Ju],x)\\
& =-\omega_2([x,y],Ju)+\omega_2(\ncp_xJu,y)-\omega_2(\ncp_yJu,x)\\
& =\omega_1([x,y],u)-\omega_1(\ncp_xu,y)+\omega_1(\ncp_yu,x)\\
&
=\omega_+([x,y],u)-\omega_+(\nabla^+_xu,y)+\omega_+(\nabla^+_yu,x)\\
&
=\omega_+([x,y],u)-\omega_+(\nabla^+_xy,u)+\omega_+(\nabla^+_yx,u)\\
& =0
\end{align*}
because of (\ref{split}), $\riii$ and since $\nabla^+$ is
torsion-free. Now suppose that $x\in\ggo_+$ and $y=Jv,\,z=Ju$ with
$u,v\in\ggo_+$. We have that \begin{align*} \dif\omega_2)(x,Jv,Ju)
& =
-\omega_2([x,Jv],Ju)+\omega_2([x,Ju],Jv)-\omega_2([Jv,Ju],x)\\
&
=\omega_2(\ncp_{Jv}x,Ju)-\omega_2(\ncp_{Ju}x,Jv)-\omega_2(J[Jv,Ju],Jx)\\
& =-\omega_1(\ncp_{Jv}x,u)+\omega_1(\ncp_{Ju}x,v)-\omega_1(J[Jv,Ju],x)\\
& =-\omega_1(\ncp_{Jv}Jx,Ju)+\omega_1(\ncp_{Ju}Jx,Jv)+\omega_1([Jv,Ju],Jx)\\
& =-\omega_1(\ncp_{Jv}Ju,Jx)+\omega_1(\ncp_{Ju}Jv,Jx)+\omega_1([Jv,Ju],Jx)\\
& =-\omega_-(\nabla^-_{Jv}Ju,Jx)+\omega_-(\nabla^-_{Ju}Jv,Jx)+\omega_1([Jv,Ju],Jx)\\
& =0 \end{align*} because of (\ref{split}), $\riii$ and since
$\nabla^-$ is torsion-free. Hence, $\riv$ holds.
\end{proof}

\medskip

As a consequence, if, for a metric $g$ on $\ggo$ compatible with
the complex product structure $\{J,E\}$, the K\"ahler form
$\omega_1$ is closed (and hence symplectic), then $\{J,E,g\}$ is a
hypersymplectic structure on $\ggo$.

\

Summing up, we have the following result, which describes the
structure of a Lie algebra admitting a hypersymplectic structure.

\begin{thm}\label{jeg}
Let $\{J,E,g\}$ be a hypersymplectic structure on $\ggo$ and let
$(\ggo,\ggo_+,\ggo_-)$ be the double Lie algebra associated to the
complex product structure $\{J,E\}$. Then $\ggo_+$ carries a flat
torsion-free connection $\nabla^+$ and a compatible symplectic
form $\omega_+$, and similarly, $\ggo_-$ carries a flat
torsion-free connection $\nabla^-$ and a compatible symplectic
form $\omega_-$. These symplectic forms are related by:
$\omega_+(x,y)=\omega_-(Jx,Jy)$ for $x,y\in\ggo_+$.
\end{thm}

\

In the next result we exhibit a converse for Theorem \ref{jeg},
showing a method to produce hypersymplectic Lie algebras beginning
with two Lie algebras equipped with compatible flat torsion-free
connections and symplectic forms. Even though the hypotheses
appearing in the statement of the theorem seem rather complicated,
we will be able to use this result to obtain all the
hypersymplectic 4-dimensional Lie algebras. This will be done in
subsequent sections.

\vs

\begin{thm}\label{const}
Consider the following data:
\begin{enumerate}
\item $\ug$ is a Lie algebra equipped with a flat torsion-free
connection $\nabla$ and a symplectic form $\omega$ such that
$\nabla\omega=0$.
\item $\vg$ is a Lie algebra equipped with a
flat torsion-free connection $\nabla'$ and a symplectic form
$\omega'$ such that $\nabla'\omega'=0$.
\item There exists a linear isomorphism $\varphi:\ug\longrightarrow\vg$ such that
\begin{enumerate}
\item[\ri] the representations $\rho:\ug\longrightarrow\glg(\vg)$ and
$\mu:\vg\longrightarrow\glg(\ug)$ defined by
\[ \rho(x)a=\varphi\nabla_x\varphi^{-1}(a),\quad \mu(a)x=\varphi^{-1}\nabla'_a\varphi(x)\]
satisfy
\begin{gather*}
\rho(x)[a,b]-[\rho(x)a,b]-[a,\rho(x)b]+\rho(\mu(a)x)b-\rho(\mu(b)x)a=0,
\label{jacobi1}\\
\mu(a)[x,y]-[\mu(a)x,y]-[x,\mu(a)y]+\mu(\rho(x)a)y-\mu(\rho(y)a)x=0,\label{jacobi2}
\end{gather*}
for all $x,y\in\ug$ and $a,b\in\vg$.
\item[\rii] $\omega(x,y)=\omega'(\varphi(x),\varphi(y))$ for all $x,y\in\ug$.
\end{enumerate}
\end{enumerate}
In this situation the vector space $\ggo=\ug\oplus\vg$ admits a
Lie bracket extending the Lie brackets on $\ug$ and $\vg$ and
there is a hypersymplectic structure on $\ggo$ such that
$\ggo_+=\ug$ and $\ggo_-=\vg$.
\end{thm}

\begin{proof}
Condition $\ri$ in the statement means that ($\ug,\vg,\rho,\mu$)
is a {\em matched pair} of Lie algebras (see \cite{Maj} or
\cite{Mas}). Thus, the bracket on $\ggo$ given by
\[ [(x,a),(y,b)]=([x,y]+\mu(a)y-\mu(b)x,[a,b]+\rho(x)b-\rho(y)a),\]
for $x,y\in\ug$ and $a,b\in\vg$ satisfies the Jacobi identity;
$\ggo$ with this Lie algebra structure will be denoted
$\ggo=\ug\bowtie\vg$ and will be called the {\em bicrossproduct}
of $\ug$ and $\vg$. Observe that $\ug$ and $\vg$ are Lie
subalgebras of $\ggo$. If we take into account the definition of
$\rho$ and $\mu$, we get that the Lie bracket between an element
of $\ug$ and one of $\vg$ is given by
\[ [(x,0),(0,a)]=(-\varphi^{-1}\nabla'_a\varphi(x),\varphi\nabla_x\varphi^{-1}(a)) \]
for $x\in\ug$ and $a\in\vg$. It has already been proved in
\cite{AS} that $\ggo=\ug\bowtie\vg$ admits a complex product
structure $\{J,E\}$, where the endomorphisms $J$ and $E$ are
defined by
\[ J(x,a)=(-\varphi^{-1}(a),\varphi(x)),\qquad E|_{\ug}=\Id,\quad E|_{\vg}=-\Id, \]
for $x\in\ug,\,a\in\vg$. Furthermore, if $\ncp$ denotes the
torsion-free connection associated to $\{J,E\}$, then the
restrictions of $\ncp$ to $\ggo_+=\ug$ and $\ggo_-=\vg$ are
precisely the original connections $\nabla$ and $\nabla'$,
respectively.

We proceed now to define a metric $g$ on $\ggo$ which will be
shown to be hypersymplectic. Let $g$ be given by
\[ g(\ug,\ug)=0,\quad g(\vg,\vg)=0, \quad g((x,a),(y,b))=\omega(\varphi^{-1}(b),x)+
\omega(\varphi^{-1}(a),y) \]
for $x,y\in\ug,\,a,b\in\vg$. It is clear that $g$ is a metric on
$\ggo$. We should check now that it satisfies (\ref{metric}). We
begin with
\begin{align*}
g(J(x,a),J(y,b)) & = g((-\varphi^{-1}(a),\varphi(x)),(-\varphi^{-1}(b),\varphi(y)) \\
                 & = -\omega(y,\varphi^{-1}(a))-\omega(x,\varphi^{-1}(b)) \\
                 & = \omega(\varphi^{-1}(a),y)+\omega(\varphi^{-1}(b),x) \\
                 & = g((x,a),(y,b))
\end{align*}
and now
\begin{align*}
g(E(x,a),E(y,b)) & = g((x,-a),(y,-b)) \\
                 & = \omega(-\varphi^{-1}(b),x)+\omega(-\varphi^{-1}(a),y)\\
                 & = -\omega(\varphi^{-1}(b),x)-\omega(\varphi^{-1}(a),y)\\
                 & = -g((x,a),(y,b)).
\end{align*}
Thus (\ref{metric}) holds and $g$ is compatible with $\{J,E\}$.

To see that with this metric we obtain a hypersymplectic structure
on $\ggo$, we only have to see that $\riii$ of Proposition
\ref{omegas} holds. Let us determine firstly the 2-form $\omega_1$
on $\ggo$:
\begin{align*}
\omega_1((x,a),(y,b)) & = g(J(x,a),(y,b))\\
& = g((-\varphi^{-1}(a),\varphi(x)),(y,b))\\
& = \omega(\varphi^{-1}(b),-\varphi^{-1}(a))+\omega(x,y)\\
& = \omega(x,y)+\omega'(a,b).
\end{align*} Therefore, the restrictions of $\omega_1$ to
$\ggo_+=\ug$ and $\ggo_-=\vg$ are precisely the original
symplectic forms $\omega$ and $\omega'$, respectively. As
$\nabla^+\omega_+=\nabla\omega=0$ and
$\nabla^-\omega_-=\nabla'\omega'=0$, we have that $g$ is a
hypersymplectic metric on $\ggo$.
\end{proof}

\

Any metric $g$ on a Lie algebra $\ggo$ determines by
left-translations a left-invariant metric on $G$, where $G$ is the
only simply connected Lie group with $\Lie(G)=\ggo$. It is easy to
verify that the Levi-Civita connection on the manifold $G$ is also
left-invariant, and hence it is determined by its values at
$\ggo\cong \tang_eG$. Therefore, the metric $g$ on $\ggo$
determines a connection $\nabla^g$ on $\ggo$, also called the
Levi-Civita connection associated to $g$. This Levi-Civita
connection is the only connection on $\ggo$ such that $\ri$ it is
torsion-free, and $\rii$ the endomorphisms
$\nabla^g_x,\,x\in\ggo$, are skew-adjoint with respect to $g$.
Just as in the positive definite case, in the neutral setting one
can prove the following equivalences:

\begin{prop}\label{equiv-hyper}
Let $\ggo$ be a Lie algebra with a complex product structure
$\{J,E\}$ and a compatible metric $g$. Let $\nabla^g$ denote the
Levi-Civita connection on $\ggo$ associated to $g$ and let
$\omega_i,\,i=1,2,3$, be the 2-forms on $\ggo$ given in
(\ref{forms}). Then the following statements are equivalent:
\begin{enumerate}
\item[\ri] The metric $g$ is hypersymplectic, i.e.,
$\dif\omega_i=0$ for $i=1,2,3$. \item[\rii] The endomorphisms $J$
and $E$ are $\nabla^g$-parallel: $\nabla^gJ=\nabla^gE=0$.
\item[\riii] The 2-forms $\omega_i,\,i=1,2,3$, are
$\nabla^g$-parallel: $\nabla^g\omega_i=0$ for $i=1,2,3$.
\end{enumerate}
\end{prop}

\vs

On any Lie algebra $\ggo$ with a hypersymplectic structure
$\{J,E,g\}$ we have canonically defined two torsion-free
connections: the connection $\ncp$ determined by the complex
product structure $\{J,E\}$ and the Levi-Civita connection
$\nabla^g$ corresponding to $g$. Recalling that $\ncp$ is the only
torsion-free connection such that $J$ and $E$ are parallel, and
taking into account the equivalence $\ri\Leftrightarrow\rii$ of
Proposition \ref{equiv-hyper}, we obtain that $\nabla^g=\ncp$.

\

We consider now the question of equivalences between
hypersymplectic structures. We have the following definition.

\begin{defi}
Let $\{J,E,g\}$ and $\{J',E',g'\}$ be hypersymplectic structures on
the Lie algebras $\ggo$ and $\ggo'$ respectively. These structures
are said to be {\em equivalent} if there exists a Lie algebra
isomorphism $\xi:\ggo\longrightarrow\ggo'$ such that
\begin{equation}\label{isom} \xi J=J'\xi,\quad \xi E=E'\xi \quad \text{and}
\quad g'(\xi x,\xi y)=g(x,y) \end{equation}
for all $x,y\in\ggo$.
\end{defi}

\medskip

\begin{rem}
The first two conditions in (\ref{isom}) mean that the underlying
complex product structures $\{J,E\}$ and $\{J',E'\}$ are
equivalent. The third condition means that $\xi$ is an isometry
between $g$ and $g'$.
\end{rem}

\

\begin{lem}
With notation as in the previous definition, let $\nabla^g$ and
$\nabla^{g'}$ be the Levi-Civita connection of $g$ and $g'$
respectively. Then $\xi$ gives an equivalence between these two
connections. Furthermore, if $\omega_i,\,i=1,2,3$, are given as in
(\ref{forms}) and $\omega'_i,\,i=1,2,3$, are defined similarly for
$\ggo'$, then $\omega_i(x,y)=\omega'_i(\xi x,\xi y)$ for all
$x,y\in\ggo$.
\end{lem}

\begin{proof}
The Levi-Civita connection $\nabla^g$ of $g$ is the only
torsion-free connection on $\ggo$ such that $\nabla^g J=\nabla^g
E=0$, and a similar statement holds for $\nabla^{g'}$. We would
like to show that $\xi \nabla^g_xy=\nabla^{g'}_{\xi x}\xi y$ for
all $x,y\in \ggo$. To see this, define a connection
$\tilde{\nabla}$ on $\ggo$ by
\[ \tilde{\nabla}_xy:=\xi^{-1}\nabla^{g'}_{\xi x}\xi y,\quad x,y\in\ggo.\]
Let us see that it is torsion-free:
\[ \tilde{\nabla}_xy-\tilde{\nabla}_yx=\xi^{-1}(\nabla^{g'}_{\xi x}\xi y-
\nabla^{g'}_{\xi y}\xi x)=\xi^{-1}[\xi x,\xi y]=[x,y]\]
for any $x,y\in\ggo$. Let us verify now that $J$ is
$\tilde{\nabla}$-parallel.
\[ \tilde{\nabla}_xJy=\xi^{-1}\nabla^{g'}_{\xi x}\xi Jy=\xi^{-1}\nabla^{g'}_{\xi x}J'\xi y=
\xi^{-1}J'\nabla^{g'}_{\xi x}\xi y=J\xi^{-1}\nabla^{g'}_{\xi x}\xi y=J\tilde{\nabla}_xy \]
for all $x,y\in\ggo$. In the same way, it can be seen that
$\tilde{\nabla}E=0$. By uniqueness, we have that
$\nabla^g=\tilde{\nabla}$ and hence $\nabla^g$ and $\nabla^{g'}$
are equivalent.

Let us check now the assertions about the symplectic forms. Let
us consider first the 2-form $\omega_1$. We have
\[\omega_1(x,y)=g(Jx,y)=g'(\xi Jx,\xi y)=g'(J'\xi x,\xi y)=\omega'_1(\xi x,\xi y)\]
for all $x,y\in\ggo$. In a similar fashion one can prove the
corresponding statements for $\omega_2$ and $\omega_3$.
\end{proof}

\vs

Motivated by the previous result, we introduce the following definition.

\begin{defi}
Let $\ggo$ be a Lie algebra equipped with a connection $\nabla$
and a symplectic form $\omega$ such that $\nabla\omega=0$, and
similarly for a Lie algebra $\ggo'$ with $\nabla'$ and $\omega'$.
We will say that $(\nabla,\omega)$ and $(\nabla,\omega)$ are {\em
symplectically equivalent} if there exists a Lie algebra
isomorphism $\xi:\ggo\longrightarrow\ggo'$ such that
\[ \xi \nabla_xy=\nabla'_{\xi x}\xi y, \qquad \omega(x,y)=\omega'(\xi x,\xi y) \]
for all $x,y\in\ggo$.
\end{defi}

\

\begin{prop}
Keep the notation from Theorem \ref{const}. Suppose that $(\nabla,
\omega)$ is symplectically equivalent to
$(\overline{\nabla},\overline{\omega})$, where $\overline
{\nabla}$ is a flat torsion-free connection on $\ug$ and
$\overline{\omega}$ is a symplectic form on $\ug$ such that
$\overline{\nabla}\overline{\omega}=0$. Equally, let
$(\nabla',\omega')$ be symplectically equivalent to
$(\overline{\nabla}',\overline{\omega}')$. Then we obtain a
matched pair of Lie algebras
$(\ug,\vg,\overline{\rho},\overline{\mu})$ and the bicrossproduct
$\overline{\ggo}=\ug
\bowtie^{\overline{\rho}}_{\overline{\mu}}\vg$ has a
hypersymplectic structure equivalent to the one on $\ggo=\ug
\bowtie^{\rho}_{\mu}\vg$.
\end{prop}

\begin{proof}
Let $\xi:\ug\longrightarrow\ug$ and $\xi':\vg\longrightarrow\vg$
be the Lie algebra isomorphisms which give the symplectic
equivalences between $(\nabla,\omega)$ and
$(\overline{\nabla},\overline{\omega})$ and between
$(\nabla',\omega')$ and $(\overline{\nabla}',\overline{\omega}')$.
Consider now the linear isomorphism $\psi:\ug\longrightarrow\vg$
given by $\psi=\xi'\varphi\xi^{-1}$. Associated to the isomorphism
$\psi$ we have the representations
$\overline{\rho}:\ug\longrightarrow\glg(\vg)$ and
$\overline{\mu}:\vg\longrightarrow\glg(\ug)$ defined by
\[ \overline{\rho}(x)a=\psi\overline{\nabla}_x\psi^{-1}(a),\quad
\overline{\mu}(a)x=\psi^{-1}\overline{\nabla}'_a\psi(x).\] It is
easily verified that $(\ug,\vg,\overline{\rho},\overline{\mu})$ is
a matched pair of Lie algebras, using that $(\ug,\vg,\rho,\mu)$ is
another matched pair of Lie algebras. We may form now the
bicrossproduct Lie algebras $\ggo=\ug \bowtie^{\rho}_{\mu}\vg$ and
$\overline{\ggo}=\ug\bowtie^{\overline{\rho}}_{\overline{\mu}}\vg$.
Furthermore, it is easy to see that
$\overline{\omega}(x,y)=\overline{\omega}'(\psi(x),\psi(y))$ for
all $x,y\in\ug$. From Theorem \ref{const}, both $\ggo$ and
$\overline{\ggo}$ have a hypersymplectic structure. Let us see now
that they are equivalent. Consider the linear isomorphism
$\eta:=\xi\oplus\xi':\ggo\longrightarrow\overline{\ggo}$ and
observe that
\begin{align*}
[\eta(x,a),\eta(y,b)] & =[(\xi x,\xi' a),(\xi y, \xi' b)] \\
& =([\xi x, \xi y]+\overline{\mu}(\xi'a)\xi y-\overline{\mu}(\xi'b)\xi x,
[\xi'a,\xi'b]+\overline{\rho}(\xi x)\xi'b-\overline{\rho}(\xi y)\xi'a)\\
& =(\xi[x,y]+\xi(\mu(a)y)-\xi(\mu(b)x),\xi'[a,b]+\xi'(\rho(x)b)-\xi'(\rho(y)a))\\
& =\eta([x,y]+\mu(a)y-\mu(b)x,[a,b]+\rho(x)b-\rho(y)a)\\
& =\eta[(x,a),(y,b)].
\end{align*}
Thus, $\eta$ is a Lie algebra isomorphism. Now,
\[ \overline{J}\eta(x,a)=\overline{J}(\xi x,\xi'a)=(-\psi^{-1}\xi'a,\psi\xi x)=
(-\xi\varphi^{-1}a,\xi'\varphi x)=\eta(-\varphi^{-1}a,\varphi x)=\eta J(x,a)\]
and
\[ \overline{E}\eta(x,a)=\overline{E}(\xi x,\xi'a)=(\xi x,-\xi'a)=\eta(x,-a)=\eta E(x,a)\]
for all $(x,a)\in\ggo$. Hence, $\eta J=\overline{J}\eta$ and $\eta
E=\overline{E}\eta$. Finally, we show that $\eta$ is an isometry
between $g$ and $\overline{g}$:
\begin{align*}
\overline{g}(\eta(x,a),\eta(y,b)) & =\overline{g}((\xi x, \xi'a),(\xi y,\xi'b))=
\omega(\psi^{-1}\xi'b, \xi x)+\omega(\psi^{-1}\xi'a,\xi y) \\
& =\omega(\xi\varphi^{-1}b,\xi x)+\omega(\xi\varphi^{-1}a,\xi y)\\
& =\omega(\varphi^{-1}b,x)+\omega(\varphi^{-1}a,y)\\
& =g((x,a),(y,b))
\end{align*}
for all $(x,a),(y,b)\in\ggo$. Thus, $\eta$ gives an equivalence
between the hypersymplectic structures on $\ggo$ and
$\overline{\ggo}$.
\end{proof}

\

\section{Symplectic flat torsion-free connections on $\RR^2$ and $\aff(\RR)$}

In the next section, we will determine all the 4-dimensional Lie
algebras which carry a hypersymplectic structure. In order to do
so, we will need to know all the flat torsion-free connections
that preserve a symplectic form on the 2-dimensional Lie algebras.
We recall that, up to isomorphism, there are only two
2-dimensional Lie algebras, namely, $\RR^2$ and the Lie algebra
$\aff(\RR)$, which has a basis $\{e_1,e_2\}$ such that
$[e_1,e_2]=e_2$. $\aff(\RR)$ is the Lie algebra of the Lie group
$\Aff(\RR)$ of affine motions of the real line.

We start with the abelian Lie algebra $\RR^2$.

\

\begin{thm}\label{abelian}
Let $\RR^2=\text{span}\{e_1,e_2\}$ denote the 2-dimensional
abelian Lie algebra and let $\omega=e^1\wedge e^2$ be the
canonical symplectic form on $\RR^2$. Then the only non zero flat
torsion-free connections $\nabla$ on $\RR^2$ such that
$\nabla\omega=0$ are the following:
\begin{enumerate}
\item[(a)] $\nabla_{e_1}e_1=\alpha e_2\;(\alpha\neq 0)$, the other possibilities being 0;
\item[(b)] $\nabla_{e_2}e_2=\alpha e_1\;(\alpha\neq 0)$, the other possibilities being 0;
\item[(c)] For $\alpha\neq 0,\beta\neq 0$,
\begin{align*}
\nabla_{e_1}e_1 &=\alpha e_1+\beta e_2,\\
\nabla_{e_1}e_2 &=-\frac{\alpha}{\beta}(\alpha e_1+\beta e_2)=\nabla_{e_2}e_1,\\
\nabla_{e_2}e_2 &=\frac{\alpha^2}{\beta^2}(\alpha e_1+\beta e_2).
\end{align*}
\end{enumerate}
\end{thm}

\begin{proof}
Let us denote
\begin{align*}
\nabla_{e_1}e_1 &=ae_1+be_2,\\
\nabla_{e_1}e_2 &=ce_1+de_2=\nabla_{e_2}e_1,\\
\nabla_{e_2}e_2 &=ge_1+he_2,
\end{align*}
with $a,b,c,d,g,h\in\RR$. Since $\nabla$ is flat, we have that
$\nabla_{e_1}\nabla_{e_2}=\nabla_{e_2}\nabla_{e_1}$, and from this
condition we obtain that
\begin{gather}
bg=cd,\nonumber\\
bc-bh+d^2-ad=0,\label{flat}\\
ag-dg+ch-c^2=0. \nonumber
\end{gather}
Now, the condition $\nabla\omega=0$ holds if and only if
$\omega(\nabla_xy,z)=\omega(\nabla_xz,y)$ for all $x,y,z\in\RR^2$.
From this we get
\[ d=-a \quad \text{and} \quad h=-c.\]
Substituting into (\ref{flat}), we obtain
\begin{gather}
bg=-ac,\nonumber\\
a^2=-bc,\label{flat2}\\
c^2=ag. \nonumber
\end{gather}

If $a=0$, then $c=0$ and $bg=0$. As $\nabla\neq 0$, then $b\neq 0$
or $g\neq 0$. If $b\neq 0$, then $g=0$ and $\nabla$ is of type (a)
in the statement. If $g\neq 0$, then $b=0$ and $\nabla$ is of type
(b) in the statement.

Let us suppose now $a\neq 0$. Then $bcg\neq 0$ and from
(\ref{flat2}) we obtain $c=-\frac{a^2}{b}$ and
$g=\frac{a^3}{b^2}$. Therefore, $\nabla$ is of type (c) in the
statement.
\end{proof}

\

In the next proposition we study the equivalences among the
connections obtained in Theorem \ref{abelian}

\begin{prop}\label{equiv-ab}
Let $\nabla$ be a non zero flat torsion-free connection on $\RR^2$ and
$\omega$ a $\nabla$-parallel symplectic form on $\RR^2$. Then
$(\nabla,\omega)$ is symplectically equivalent to $(\nabla^0,
e^1\wedge e^2)$, where $\{e_1,e_2\}$ is a suitable basis of
$\RR^2$, $\{e^1,e^2\}$ is the dual basis and $\nabla^0$ is given
by
\[ \nabla^0_{e_1}e_1=e_2,\,\nabla^0_{e_1}e_2=0,\,\nabla^0_{e_2}\equiv 0.\]
This flat torsion-free connection on $\RR^2$ is complete.
\end{prop}

\begin{proof}
There exists a basis $\{e_1,e_2\}$ of $\RR^2$ such that
$\omega=e^1\wedge e^2$. Since $\nabla\omega=0$, the connection
$\nabla$ must be one of those given by Theorem \ref{abelian}.

Let us suppose first that $\nabla$ is of type (a) in Theorem
\ref{abelian}. The linear isomorphism of $\RR^2$ which gives the
symplectic equivalence between $\nabla$ and $\nabla^0$ is given by
\[ \xi=\begin{pmatrix} \alpha^{1/3} & 0 \cr 0 & \alpha^{-1/3}\cr \end{pmatrix}\]
in the ordered basis $\{e_1,\,e_2\}$.

Suppose now that the connection $\nabla$ is of type (b) in Theorem
\ref{abelian}. The linear isomorphism of $\RR^2$ which gives the
symplectic equivalence between $\nabla$ and $\nabla^0$ is given by
\[ \xi=\begin{pmatrix} 0 & -\alpha^{1/3} \cr \alpha^{-1/3} & 0\cr \end{pmatrix}\]
in the ordered basis $\{e_1,\,e_2\}$.

Finally, if $\nabla$ is of type (c) in Theorem \ref{abelian}, we
may take the following isomorphism of $\RR^2$:
\[ \xi=\begin{pmatrix} \beta^{1/3} & -\alpha\beta^{-2/3} \cr 0 & \beta^{-1/3}\cr
\end{pmatrix}.\]
The verification of all these statements is a simple matter. We
would like now to check that the connection $\nabla^0$ is
complete. In order to do so, we will use equation
(\ref{complete}). Let $x(t)=a_1(t)e_1+a_2(t)e_2$ be a curve on
$\ggo$ which satisfies $\dot{x}(t)=-\nabla^0_{x(t)}x(t)$. Thus, we
obtain the system of differential equations
\[ \begin{cases}\dot{a_1}=0,\\
\dot{a_2}=-a_1^2.
\end{cases} \]
The solutions of this system are clearly defined for every
$t\in\RR$ and therefore $\nabla^0$ is complete.
\end{proof}

\vs

\begin{rem}
In \cite{RG}, a classification of flat torsion-free connections
(up to equivalence) on $\RR^2$ is given. The flat torsion-free
connection $\nabla^0$ from Proposition \ref{equiv-ab} belongs to
the class $A_4$ of that classification.
\end{rem}

\

Next, we move on to consider the other 2-dimensional Lie algebra,
$\aff(\RR)$.

\

\begin{thm}\label{aff}
Let $\aff(\RR)=\text{span}\{e_1,e_2\}$ denote the 2-dimensional
Lie algebra with Lie bracket $[e_1,e_2]=e_2$ and let $\omega=e^1\wedge
e^2$ be the canonical symplectic form on $\aff(\RR)$. Then the
only flat torsion-free connections $\nabla$ on
$\aff(\RR)$ such that $\nabla\omega=0$ are the following:
\begin{enumerate}
\item[(a)] For $\alpha\in\RR$,
\begin{align*}
\nabla_{e_1}e_1 &=-e_1+\alpha e_2, \\
\nabla_{e_1}e_2 &=e_2, \\
\nabla_{e_2}e_1 &=\nabla_{e_2}e_2=0.
\end{align*}
\item[(b)] For $\alpha\in\RR$,
\begin{align*}
\nabla_{e_1}e_1 &=-\frac{1}{2}e_1+\alpha e_2, \\
\nabla_{e_1}e_2 &=\frac{1}{2}e_2, \\
\nabla_{e_2}e_1 &=-\frac{1}{2}e_2,\\
\nabla_{e_2}e_2 &=0.
\end{align*}
\end{enumerate}
\end{thm}

\begin{proof}
Let us denote
\begin{align*}
\nabla_{e_1}e_1 &=ae_1+be_2,\\
\nabla_{e_1}e_2 &=ce_1+de_2,\\
\nabla_{e_2}e_2 &=ge_1+he_2,
\end{align*}
with $a,b,c,d,g,h\in\RR$. Since $\nabla$ is torsion-free, we have
\[\nabla_{e_2}e_1=ce_1+(d-1)e_2.\]

The condition $\nabla\omega=0$ implies that $d=-a$ and $h=-c$.
Taking this into account and using that $\nabla$ is flat, we
obtain the following equations
\begin{gather}
c(a+2)+bg=0, \nonumber\\
g(2a-1)-2c^2=0, \label{afin}\\
2bc+(a+1)(2a+1)=0. \nonumber
\end{gather}
From the third equation in (\ref{afin}) we get
\begin{equation}\label{bc}
2a^2+3a+(2bc+1)=0.
\end{equation}
Also, from (\ref{afin}) we see immediately that $a\neq
\frac{1}{2}$. Hence $g=\frac{2c^2}{2a-1}$ and substituting into
the first equation we have
\[ c\left((a+2)+\frac{2bc}{2a-1}\right)=0.\]
If $c\neq 0$, then $(a+2)(2a-1)+2bc=0$ and $2a^2+3a+2bc-2=0$,
which combined with (\ref{bc}) yields a contradiction. Thus, $c=0$
and the system (\ref{afin}) becomes
\begin{gather}
bg=0, \nonumber\\
g(2a-1)=0, \\
(a+1)(2a+1)=0. \nonumber
\end{gather}
Therefore, $g=0$ (since $a\neq \frac{1}{2}$), $b\in\RR$ is
arbitrary and $a=-1$ or $a=-\frac{1}{2}$. In the first case, we
obtain a connection of type (a) and in the second case we obtain a
connection of type (b). The proof is complete.
\end{proof}

\

In the next proposition we deal with the equivalences of the
connections obtained in Theorem \ref{aff}.

\

\begin{prop}\label{equiv-aff}
Let $\nabla$ be a flat torsion-free connection on $\aff(\RR)$ and
$\omega$ a $\nabla$-parallel symplectic form on $\aff(\RR)$. Then
$(\nabla,\omega)$ is symplectically equivalent to either
$(\nabla^1, e^1\wedge e^2)$ or $(\nabla^2, e^1\wedge e^2)$, where
$\{e_1,e_2\}$ is a suitable basis of $\aff(\RR)$, $\{e^1,e^2\}$ is
the dual basis and $\nabla^1,\,\nabla^2$ are given by:
\[ \nabla^1_{e_1}e_1=-e_1,\quad \nabla^1_{e_1}e_2=e_2,\quad \nabla^1_{e_2}\equiv 0 \]
and
\[ \nabla^2_{e_1}e_1=-\frac{1}{2}e_1,\quad \nabla^2_{e_1}e_2=\frac{1}{2}e_2,\quad
\nabla^2_{e_2}e_1=-\frac{1}{2}e_2, \quad \nabla^2_{e_2}e_2=0.\]
Both connections $\nabla^1$ and $\nabla^2$ on $\aff(\RR)$ are not complete.
\end{prop}

\begin{proof}
Let $\{{\tilde e}_1,{\tilde e}_2\}$ be a basis of $\aff(\RR)$ such
that $[{\tilde e}_1,{\tilde e}_2]={\tilde e}_2$. There exists
$\lambda\neq 0$ such that $\omega=\lambda({\tilde e}^1\wedge
{\tilde e}^2)$. Set $e_1:={\tilde e}_1,\,e_2:=\lambda{\tilde
e}_2$. We have then $[e_1,e_2]=e_2$ and
\[ \omega=\lambda({\tilde e}^1\wedge {\tilde e}^2)=\lambda(e^1\wedge \lambda^{-1}e^2)=
e^1\wedge e^2.\]
So, we have $\nabla(e^1\wedge e^2)=0$, and then $\nabla$ must be
one of the flat torsion-free connections given in Theorem
\ref{aff}.

Let us suppose first that $\nabla$ is of type (a) in Theorem
\ref{aff}. The linear isomorphism of $\aff(\RR)$ which gives the
symplectic equivalence between $\nabla$ and $\nabla^1$ is given by
\[ \xi=\begin{pmatrix} 1 & 0 \cr \frac{1}{2}\alpha & 1 \cr \end{pmatrix}\]
in the ordered basis $\{e_1,\,e_2\}$.

If we take now a connection $\nabla$ of type (b) in Theorem
\ref{aff}, the linear isomorphism of $\aff(\RR)$ giving the
symplectic equivalence between $\nabla$ and $\nabla^2$ is
\[ \xi=\begin{pmatrix} 1 & 0 \cr 2\alpha & 1 \cr \end{pmatrix}\]
in the ordered basis $\{e_1,\,e_2\}$.

Next, we observe that $\nabla^1$ and $\nabla^2$ are not
equivalent. If they were, the subspaces
$W_1=\{x\in\aff(\RR):\nabla^1_x\equiv 0\}$ and
$W_2=\{x\in\aff(\RR):\nabla^2_x\equiv 0\}$ of $\aff(\RR)$ should
be isomorphic. However, it is clear that $\dim W_1=1$ while
$W_2=\{0\}$. Thus, these two connections are not equivalent.

Finally, we show that these connections are not complete. Suppose
$x(t)=a_1(t)e_1+a_2(t)e_2$ is a curve on $\aff(\RR)$ that
satisfies $\dot{x}(t)=-\nabla^1_{x(t)}x(t)$. Thus, we obtain the
system of differnetial equations
\[ \begin{cases}\dot{a_1}=a_1^2,\\
\dot{a_2}=-a_1a_2.
\end{cases} \]
From the first equation in the system we obtain that $a_1(t)$
cannot be defined in the whole real line; thus $\nabla^1$ is not
complete. Analogously, if $x(t)=a_1(t)e_1+a_2(t)e_2$ is a curve on
$\aff(\RR)$ that satisfies $\dot{x}(t)=-\nabla^2_{x(t)}x(t)$, we
have the system
\[ \begin{cases}\dot{a_1}=\ft a_1^2,\\
\dot{a_2}=0.
\end{cases} \]
We obtain again that $a_1(t)$ cannot be defined in the whole real
line; thus $\nabla^2$ is not complete.
\end{proof}

\

\section{Hypersymplectic 4-dimensional Lie algebras}

In this section we will determine all $4$-dimensional Lie algebras
which carry a hypersymplectic structure, by employing Theorem
\ref{const}. We will also be able to obtain a parametrization of
the underlying complex product structures, up to equivalence. In
the next section we will exhibit explicit descriptions of the
hypersymplectic metrics in each case.

A classification of 4-dimensional Lie algebras admitting a complex product structure was given by Blazi\'c and Vukmirovi\'c in \cite{BV}, where they refer to complex product structures as para-hypercomplex structures. The family of Lie algebras that we will obtain below can be found within this classification.

We introduce first some notation. We will consider the following
$4$-dimensional Lie algebras:

\n $\bullet\,\ggo^h_0=\text{span}\{v_0,v_1,v_2,v_3\}$ with
$[v_1,v_2]=v_3$,

\n $\bullet\,\ggo^h_1=\text{span}\{v_0,v_1,v_2,v_3\}$ with
$[v_0,v_1]=v_1,\, [v_0,v_2]=-v_2,\,[v_0,v_3]=-v_3$,

\n $\bullet\,\ggo^h_2=\text{span}\{v_0,v_1,v_2,v_3\}$ with
$[v_0,v_1]=2v_1,\,
[v_0,v_2]=-v_2,\,[v_0,v_3]=v_3,\,[v_1,v_2]=v_3$.

\begin{rems}
$\ri$ The Lie algebra $\ggo^h_0\cong\hg_3\times\RR$ is a central
extension of the $3$-dimensional Heisenberg algebra. It is the
only $2$-step nilpotent $4$-dimensional Lie algebra.

\smallbreak\noi $\rii$ The Lie algebra $\ggo^h_1$ is an extension
of $\RR^3$ and it lies in the class $\rg_{4,-1,-1}$ of the
classification of $4$-dimensional solvable Lie algebras given in
\cite{ABDO}. It is $2$-step solvable and not unimodular.

\smallbreak\noi $\riii$ The Lie algebra $\ggo^h_2$ is an extension
of $\hg_3$ and it lies in the class $\dg_{4,2}$ of the
classification of $4$-dimensional solvable Lie algebras given in
\cite{ABDO}. It is $3$-step solvable and not unimodular.
\end{rems}

\

\begin{thm}\label{four}
Let $\ggo$ be a $4$-dimensional Lie algebra carrying a
hypersymplectic structure. Then $\ggo$ is isomorphic to either
$\RR^4,\,\ggo^h_0,\,\ggo^h_1$ or $\ggo^h_2$. Furthermore, the
parametrization of the complex product structures in each case is
given by:

\smallskip

\noi $\ri$ If $\ggo\cong\RR^4$, the underlying complex product
structure is equivalent to $\{J,E\}$, where
\[ J=\begin{pmatrix} 0 & -\Id \cr \Id & 0 \cr \end{pmatrix},\quad
E=\begin{pmatrix} \Id & 0 \cr 0 & -\Id \cr \end{pmatrix}\]
(with $\Id$ the $(2\times 2)$-identity matrix) in some ordered basis of $\RR^4$.

\smallskip

\noi $\rii$ If $\ggo\cong\ggo^h_0$, then the underlying complex
product structure on $\ggo$ is equivalent to one and only one of
$\{J^{(0)},\,E^{(0)}_{\theta}\}$, where
\[ {\small J^{(0)}= \begin{pmatrix} 0 & -1 & 0 & 0 \cr
                   1 & 0 & 0 & 0 \cr
                   0 & 0 & 0 & -1 \cr
                   0 & 0 & 1 & 0 \cr
\end{pmatrix},\quad E^{(0)}_{\theta}=\begin{pmatrix}
                      1 & 0 & 0 & 0 \cr 0 & -1 & 0 & 0 \cr 0 & 0 & \cos\theta &
                      \sin\theta \cr 0 & 0 & \sin\theta & -\cos\theta \cr
                      \end{pmatrix}}\]
for $\theta\in[0,2\pi)$, in the ordered basis $\{v_1,v_2,v_3,v_0\}$.

\noi $\riii$ If $\ggo\cong\ggo^h_1$, then the underlying complex
product structure on $\ggo$ is equivalent to one and only one of
$\{J^{(1)},\,E^{(1)}_{\theta,d}\}$ or $\{J^{(1)},E^{(1)}_1\}$,
where ${\small J^{(1)}=\begin{pmatrix} 0 & -1 & 0 & 0 \cr
                   1 & 0 & 0 & 0 \cr
                   0 & 0 & 0 & -1 \cr
                   0 & 0 & 1 & 0 \cr
\end{pmatrix}}$ and
\[{\small E^{(1)}_{\theta,d}=\begin{pmatrix} \cos\theta & \sin\theta & 0 & 0 \cr
                   \sin\theta & -\cos\theta & 0 & 0 \cr
                   -d\sin\theta & d(1+\cos\theta) & 1 & 0 \cr
                    d(1+\cos\theta) & d\sin\theta & 0 & -1 \cr
                   \end{pmatrix} , \quad
E^{(1)}_1=\begin{pmatrix} -1 & 0 & 0 & 0 \cr
                 0 & 1 & 0 & 0 \cr
                 -2 & 0 & 1 & 0 \cr
                 0 & 2 & 0 & -1 \cr
\end{pmatrix}}\]
for $\theta\in[0,2\pi)$ and $d=0$ or $d=1$, in the ordered basis $\{v_0,v_1,v_2,v_3\}$.

\smallskip

\noi $\riv$ If $\ggo\cong\ggo^h_2$, then the underlying complex
product structure on $\ggo$ is equivalent to one and only one of
$\{J^{(2)},\,E^{(2)}_{\theta,d}\}$ or $\{J^{(2)},E^{(2)}_1\}$,
where ${\small J^{(2)}=\begin{pmatrix} 0 & 1 & 0 & 0 \cr
                   -1 & 0 & 0 & 0 \cr
                   0 & 0 & 0 & -1 \cr
                   0 & 0 & 1 & 0 \cr
\end{pmatrix}}$ and
\[ {\small E^{(2)}_{\theta,d}=\begin{pmatrix}   \cos\theta & \sin\theta & 0 & 0 \cr
                      \sin\theta & -\cos\theta & 0 & 0 \cr
                      d\sin\theta(1+\cos\theta) & -d\cos\theta(1+\cos\theta) &
                      \cos\theta & -\sin\theta \cr
                      d\cos\theta(1+\cos\theta) & d\sin\theta(1+\cos\theta) &
                      -\sin\theta & -\cos\theta \cr
                     \end{pmatrix},\,
E^{(2)}_1=\begin{pmatrix} -1 & 0 & 0 & 0 \cr
                 0 & 1 & 0 & 0 \cr
                 0 & 2 & -1 & 0 \cr
                 -2 & 0 & 0 & 1 \cr
\end{pmatrix}}\]
for $\theta\in[0,2\pi)$ and $d=0$ or $d=1$, in the ordered basis $\{v_0,v_2,v_1,v_3\}$.
\end{thm}

\smallskip

\begin{rem}
Note that $E^{(1)}_{\pi,0}=E^{(1)}_{\pi,1}$ and
$E^{(2)}_{\pi,0}=E^{(2)}_{\pi,1}$.
\end{rem}

\smallskip

\begin{proof}
We will construct explicitly all 4-dimensional Lie algebras
carrying a hypersymplectic structure using Theorem \ref{const}. In
order to do so, we have to determine the two 3-tuples
$(\ug,\nabla,\omega),\,(\vg,\nabla',\omega')$ and the linear
isomorphism $\varphi:\ug\longrightarrow\vg$ which satisfy the
conditions of this theorem. The Lie algebras $\ug$ and $\vg$ are
2-dimensional, and therefore they are isomorphic either to $\RR^2$
or $\aff(\RR)$. The flat torsion-free connections on these Lie
algebras which are compatible with the canonical symplectic forms
were determined in \S 4. We only have to establish the linear
isomorphisms $\varphi$ which are admissible. We will do this in
several steps.

\

{\bf Case (A): $\ug=\RR^2$ and $\vg=\RR^2$.}

We fix a basis $\{e_1,e_2\}$ of $\ug$ and its associated
symplectic form $\omega=e^1\wedge e^2$, where $\{e^1,e^2\}$ is the
dual basis. In the same way we fix a basis $\{f_1,f_2\}$ of $\vg$
and its associated symplectic form $\omega'=f^1\wedge f^2$, where
$\{f^1,f^2\}$ is the dual basis. In this case there are only two
connections to be considered: the connection identically zero and
the connection $\nabla^0$ which appears in Proposition
\ref{equiv-ab}.

\

(A1) $\nabla=0$ and $\nabla'=0$.

Here $\ggo=\ug\bowtie\vg=\RR^4$ is the abelian 4-dimensional Lie
algebra and the complex product structure is the canonical one, given
by the matrices in the statement of the theorem.

\

(A2) $\nabla=\nabla^0$ and $\nabla'=0$.

Here $\ggo=\ug\bowtie\vg=\RR^2\ltimes\RR^2$. In this special case
we may simply suppose that the linear isomorphism
$\varphi:\aff(\RR)\longrightarrow\RR^2$ that we are seeking
satisfies $\varphi(e_i)=f_i,\,i=1,2$. It is easy to see that this
isomorphism is compatible with $\nabla$ and $\nabla'$ and also
with $\omega$ and $\omega'$. Hence, we obtain a hypersymplectic
structure on $\ggo$. Let us identify this Lie algebra. If we
denote $e_i:=(e_i,0)$ and $f_i:=(0,f_i)$ for $i=1,2$, then the
only non-zero bracket is $[e_1,f_1]=f_2$. We also have $Je_i=f_i$
and $Ee_i=e_i,\,Ef_i=-f_i'$. We will make now a change of basis,
setting:
\[ v_1:=e_1,\quad v_2:=f_1,\quad v_3:=f_2,\quad v_0:=-e_2.\]
Then we have $[v_1,v_2]=v_3$ and therefore $\ggo\cong\ggo^h_0$.
The complex structure $J$ is given by $Jv_1=v_2,\,Jv_3=v_0$ and
the eigenspaces corresponding to $E$ are
$\ggo_+=\text{span}\{v_0,v_1\}$ and
$\ggo_-=\text{span}\{v_2,v_3\}$. This complex product structure is
equivalent to $\{J^{(0)},E^{(0)}_{\theta}\}$ with $\theta=\pi$.

\

(A3) $\nabla=0$ and $\nabla'=\nabla^0$.

Here $\ggo=\ug\bowtie\vg=\RR^2\ltimes\RR^2$. In this special case
we may simply suppose that the linear isomorphism
$\varphi:\aff(\RR)\longrightarrow\RR^2$ that we are seeking
satisfies $\varphi(e_i)=f_i,\,i=1,2$. It is easy to see that this
isomorphism is compatible with $\nabla$ and $\nabla'$ and also
with $\omega$ and $\omega'$. Hence, we obtain a hypersymplectic
structure on $\ggo$. Let us identify this Lie algebra. If we
denote $e_i:=(e_i,0)$ and $f_i:=(0,f_i)$ for $i=1,2$, then the
only non zero bracket is  $[e_1,f_1]=e_2$. We have that $Je_i=f_i$
and $Ee_i=e_i,\,Ef_i=-f_i$. We will make now a change of basis,
setting:
\[ v_1:=e_1,\quad v_2:=f_1,\quad v_3:=e_2,\quad v_0:=f_2.\]
Then we have $[v_1,v_2]=v_3$ and therefore $\ggo\cong\ggo^h_0$.
The complex structure $J$ is given by $Jv_1=v_2,\,Jv_3=v_0$ and
the eigenspaces corresponding to $E$ are
$\ggo_+=\text{span}\{v_1,v_3\}$ and
$\ggo_-=\text{span}\{v_0,v_2\}$. This complex product structure is
equivalent to $\{J^{(0)},E^{(0)}_{\theta}\}$ with $\theta=0$.

\

(A4) $\nabla=\nabla^0$ and $\nabla'=\nabla^0$.

We are looking for a linear isomorphism
$\varphi:\ug\longrightarrow\vg$ compatible with $\nabla$ and
$\nabla'$ and also with $\omega$ and $\omega'$. After lengthy
computations, we obtain that $\varphi$ must be of the form
$\varphi=\begin{pmatrix} a & 0 \cr b & d \cr \end{pmatrix}$, with
$ad=1$ (in the ordered bases $\{e_1,e_2\},\,\{f_1,f_2\}$). Hence,
we have a hypersymplectic structure on the bicrossproduct Lie
algebra $\ggo:=\ug\bowtie\vg=\RR^2\bowtie\RR^2$. Let us describe
this Lie algebra. Let us denote $e_i:=(e_i,0),\,f_i:=(0,f_i)$ for
$i=1,2$; the only non zero bracket is $[e_1,f_1]=-a^2e_2+d^2f_2$.
The complex product structure on this Lie algebra is given by
\[ Je_1=af_1+bf_2,\quad Je_2=df_2,\quad Jf_1=-de_1+be_2,\quad Jf_2=-ae_2 \]
and $Ee_i=e_i,\,Ef_i=-f_i$ for $i=1,2$. We will make now a change of basis, setting:
\[ v_1:=e_1,\quad v_2:=af_1+bf_2,\quad v_3:=a(-a^2e_2+d^2f_2),\quad v_0:=-a(de_2+af_2).\]
Then we have $[v_1,v_2]=v_3$ and hence $\ggo\cong\ggo^h_0$. The
complex structure $J$ is given by $Jv_1=v_2,\,Jv_3=v_0$ and the
eigenspaces corresponding to $E$ are
\[ \ggo_+=\text{span}\left\{v_1,\frac{a^3}{a^6+1}v_3+\frac{1}{a^6+1}v_0\right\},\quad
\ggo_-=\text{span}\left\{v_2,-\frac{1}{a^6+1}v_3+\frac{a^3}{a^6+1}v_0\right\}. \]
This complex product structure is equivalent to
$\{J^{(0)},E^{(0)}_{\theta}\}$, where $\theta$ is given by
$\cos(\theta/2)=\frac{a^3}{\sqrt{a^6+1}},\,\sin(\theta/2)=\frac{1}{\sqrt{a^6+1}}$.
Note that $\theta\neq 0$ and $\theta\neq\pi$.

\

{\bf Case (B): $\ug=\aff(\RR)$ and $\vg=\RR^2$.}

We fix a basis $\{e_1,e_2\}$ of $\ug$ such that $[e_1,e_2]=e_2$
and its associated symplectic form $\omega=e^1\wedge e^2$, where
$\{e^1,e^2\}$ is the dual basis. In the same way we fix a basis
$\{f_1,f_2\}$ of $\vg$ and its associated symplectic form
$\omega'=f^1\wedge f^2$, where $\{f^1,f^2\}$ is the dual basis.

\

(B1) $\nabla=\nabla^1$ and $\nabla'=0$.

In this case, we have $\ggo:=\ug\bowtie\vg=\aff(\RR)\ltimes\RR^2$.
Let us describe this Lie algebra. If we denote
$e_i:=(e_i,0),\,f_i:=(0,f_i)$ for $i=1,2$, then the only non zero
brackets are
\[ [e_1,e_2]=e_2,\quad [e_1,f_1]=-f_1,\quad [e_1,f_2]=f_2.\]
The complex product structure on this Lie algebra is given by
$Je_i=f_i$ and $Ee_i=e_i,\,Ef_i=-f_i$ for $i=1,2$. We will make a
change of basis, setting
$v_0=-e_1,\,v_1=-f_1,\,v_2=e_2,\,v_3=f_2$. Thus,
\[ [v_0,v_1]=v_1,\quad [v_0,v_2]=-v_2,\quad [v_0,v_3]=-v_3,\]
\[ Jv_0=v_1,\quad Jv_2=v_3,\]
and then $\ggo\cong\ggo^h_1$. The eigenspaces corresponding to $E$ are
\[ \ggo_+=\text{span}\{v_0,\,v_2\},\quad\ggo_-=\text{span}\{v_1,\,v_3\}.\]
This complex product structure is equivalent to
$\{J^{(1)},E^{(1)}_{\theta,0}\}$ with $\theta=0$.

\

(B2) $\nabla=\nabla^1$ and $\nabla'=\nabla^0$.

We are seeking a linear isomorphism
$\varphi:\aff(\RR)\longrightarrow\RR^2$ compatible with $\nabla$
and $\nabla'$ and also with $\omega$ and $\omega'$. It can be seen
that $\varphi$ must be of the form $\varphi=\begin{pmatrix} a & 0
\cr b & d \cr \end{pmatrix}$, with $ad=1$ (in the ordered bases
$\{e_1,e_2\},\,\{f_1,f_2\}$). Hence, we have a hypersymplectic
structure on the bicrossproduct Lie algebra
$\ggo:=\ug\bowtie\vg=\aff(\RR)\bowtie\RR^2$. Let us denote
$e_i:=(e_i,0),\,f_i:=(0,f_i)$ for $i=1,2$; the only non zero
brackets are
\[ [e_1,e_2]=e_2,\quad [e_1,f_1]=-a^2e_2-f_1-2bdf_2,\quad [e_1,f_2]=f_2.\]
The complex product structure on this Lie algebra is given by
\[ Je_1=af_1+bf_2,\quad Je_2=df_2,\quad Jf_1=-de_1+be_2,\quad Jf_2=-ae_2 \]
and $Ee_i=e_i,\,Ef_i=-f_i$ for $i=1,2$. We will make now a change
of basis, setting $v_0:=-e_1+\frac{a^2}{2}f_2,\quad
v_1:=-\frac{a^3}{2}e_2-af_1-bf_2,\quad v_2:=e_2,\quad v_3:=df_2$.
Thus,
\[ [v_0,v_1]=v_1,\quad [v_0,v_2]=-v_2,\quad [v_0,v_3]=-v_3,\]
\[ Jv_0=v_1,\quad Jv_2=v_3,\]
and then $\ggo\cong\ggo^h_1$. The eigenspaces corresponding to
$E$, in this new basis, are given by
\[ \ggo_+=\text{span}\left\{v_0-\frac{a^3}{2}v_3,\,v_2\right\},\quad
\ggo_-=\text{span}\left\{v_1+\frac{a^3}{2}v_2,\,v_3\right\}.\]
This complex product structure is equivalent to
$\{J^{(1)},E^{(1)}_{\theta,1}\}$ with $\theta=0$.

\

(B3) $\nabla=\nabla^2$ and $\nabla'=0$.

In this case, we have $\ggo:=\ug\bowtie\vg=\aff(\RR)\ltimes\RR^2$.
Let us describe this Lie algebra. If we denote
$e_i:=(e_i,0),\,f_i:=(0,f_i)$ for $i=1,2$, then the only non zero
brackets are
\[ [e_1,e_2]=e_2,\quad [e_1,f_1]=-\frac{1}{2}f_1,\quad [e_1,f_2]=\frac{1}{2}f_2,\quad
[e_2,f_1]=-\frac{1}{2}f_2.\]
The complex product structure on this Lie algebra is given by
$Je_i=f_i$ and $Ee_i=e_i,\,Ef_i=-f_i$ for $i=1,2$. We will make a
change of basis, setting
$v_0:=2e_1,\,v_1:=e_2,\,v_2:=-2f_1,\,v_3:=f_2$. Thus,
\[ [v_0,v_1]=2v_1,\quad [v_0,v_2]=-v_2,\quad [v_0,v_3]=v_3,\quad [v_1,v_2]=v_3, \]
\[ Jv_0=-v_2,\quad Jv_1=v_3, \]
and then $\ggo\cong\ggo^h_2$. The eigenspaces corresponding to
$E$, in this new basis, are given by
\[ \ggo_+=\text{span}\{v_0,\,v_1\},\quad\ggo_-=\text{span}\{v_2,\,v_3\}.\]
This complex product structure is equivalent to
$\{J^{(2)},E^{(2)}_{\theta,0}\}$ with $\theta=0$.

\

(B4) $\nabla=\nabla^2$ and $\nabla'=\nabla^0$.

We are looking for a linear isomorphism
$\varphi:\aff(\RR)\longrightarrow\RR^2$ compatible with $\nabla$
and $\nabla'$ and also with $\omega$ and $\omega'$. Again,
$\varphi$ must be of the form $\varphi=\begin{pmatrix} a & 0 \cr b
& d \cr \end{pmatrix}$, with $ad=1$ (in the ordered bases
$\{e_1,e_2\},\,\{f_1,f_2\}$). Hence, we have a hypersymplectic
structure on the bicrossproduct Lie algebra
$\ggo:=\ug\bowtie\vg=\aff(\RR)\bowtie\RR^2$. Let us denote
$e_i:=(e_i,0),\,f_i:=(0,f_i)$ for $i=1,2$; the only non zero
brackets are
\[ [e_1,e_2]=e_2,\quad [e_1,f_1]=-a^2e_2-\frac{1}{2}f_1-bdf_2,\quad
[e_1,f_2]=\frac{1}{2}f_2,\quad [e_2,f_1]=-\frac{1}{2}d^2f_2. \]
The complex product structure on this Lie algebra is given by
\[ Je_1=af_1+bf_2,\quad Je_2=df_2,\quad Jf_1=-de_1+be_2,\quad Jf_2=-ae_2 \]
and $Ee_i=e_i,\,Ef_i=-f_i$ for $i=1,2$. We will make now a change
of basis, setting $v_0:=2e_1-\frac{4}{3}a^2f_2,\quad v_1:=
e_2,\quad v_2:=-\left(\frac{4}{3}a^3e_2+2af_1+2bf_2\right),\quad
v_3:=df_2$. Thus,
\[ [v_0,v_1]=2v_1,\quad [v_0,v_2]=-v_2,\quad [v_0,v_3]=v_3,\quad [v_1,v_2]=v_3, \]
\[ Jv_0=-v_2,\quad Jv_1=v_3, \]
and then $\ggo\cong\ggo^h_2$. The eigenspaces corresponding to
$E$, in this new basis, are given by
\[ \ggo_+=\text{span}\left\{v_0+\frac{4}{3}a^3v_3,\,v_1\right\},\quad
\ggo_-=\text{span}\left\{v_2+\frac{4}{3}a^3v_1,\,v_3\right\}.\]
This complex product structure is equivalent to
$\{J^{(2)},E^{(2)}_{\theta,1}\}$ with $\theta=0$.

\

{\bf Case (B'): $\ug=\RR^2$ and $\vg=\aff(\RR)$.}

We fix a basis $\{e_1,e_2\}$ of $\ug$ and its associated
symplectic form $\omega=e^1\wedge e^2$, where $\{e^1,e^2\}$ is the
dual basis. In the same way we fix a basis $\{f_1,f_2\}$ of $\vg$
such that $[f_1,f_2]=f_2$ and its associated symplectic form
$\omega'=f^1\wedge f^2$, where $\{f^1,f^2\}$ is the dual basis.

\

(B1') $\nabla=0$ and $\nabla'=\nabla^1$.

In this case, we have $\ggo:=\ug\bowtie\vg=\aff(\RR)\ltimes\RR^2$.
If we denote $e_i:=(e_i,0),\,f_i:=(0,f_i)$ for $i=1,2$, then the
only non zero brackets are
\[ [e_1,f_1]=e_1,\quad [e_2,f_1]=-e_2,\quad [f_1,f_2]=f_2.\]
The complex product structure on this Lie algebra is given by
$Je_i=f_i$ and $Ee_i=e_i,\,Ef_i=-f_i$ for $i=1,2$. We will make a
change of basis, setting $v_0=-f_1,\,v_1=e_1,\,v_2=e_2,\,v_3=f_2$.
Thus,
\[ [v_0,v_1]=v_1,\quad [v_0,v_2]=-v_2,\quad [v_0,v_3]=-v_3,\]
\[ Jv_0=v_1,\quad Jv_2=v_3,\]
and then $\ggo\cong\ggo^h_1$. The eigenspaces corresponding to $E$
are
\[ \ggo_+=\text{span}\{v_1,\,v_2\},\quad \ggo_-=\text{span}\{v_0,\,v_3\}.\]
This complex product structure is equivalent to
$\{J^{(1)},E^{(1)}_{\theta,0}=E^{(1)}_{\theta,1}\}$ with
$\theta=\pi$.

\

(B2') $\nabla=\nabla^0$ and $\nabla'=\nabla^1$.

We are seeking a linear isomorphism
$\varphi:\aff(\RR)\longrightarrow\RR^2$ compatible with $\nabla$
and $\nabla'$ and also with $\omega$ and $\omega'$. It can be seen
that $\varphi$ must be of the form $\varphi=\begin{pmatrix} a & 0
\cr b & d \cr \end{pmatrix}$, with $ad=1$ (in the ordered bases
$\{e_1,e_2\},\,\{f_1,f_2\}$). Hence, we have a hypersymplectic
structure on the bicrossproduct Lie algebra
$\ggo:=\ug\bowtie\vg=\RR^2\bowtie\aff(\RR)$. Let us denote
$e_i:=(e_i,0),\,f_i:=(0,f_i)$ for $i=1,2$; the only non zero
brackets are
\[ [e_1,f_1]=e_1-2abe_2+d^2f_2,\quad [e_2,f_1]=-e_2,\quad [f_1,f_2]=f_2. \]
The complex product structure on this Lie algebra is given by
\[ Je_1=af_1+bf_2,\quad Je_2=df_2,\quad Jf_1=-de_1+be_2,\quad Jf_2=-ae_2 \]
and $Ee_i=e_i,\,Ef_i=-f_i$ for $i=1,2$. We will make now a change
of basis, setting
$v_0:=\frac{d^2}{2}e_2-f_1,\,v_1:=de_1-be_2+\frac{d^3}{2}f_2,\,v_2:=e_2,\,v_3:=df_2$.
Thus,
\[ [v_0,v_1]=v_1,\quad [v_0,v_2]=-v_2,\quad [v_0,v_3]=-v_3,\]
\[ Jv_0=v_1,\quad Jv_2=v_3,\]
and then $\ggo\cong\ggo^h_1$. The eigenspaces corresponding to
$E$, in this new basis, are given by
\[ \ggo_+=\text{span}\left\{v_1-\frac{d^2}{2}v_3,\,v_2\right\},\quad\ggo_-=
\text{span}\left\{-v_0+\frac{d^2}{2}v_2,\,v_3\right\}.\] This
complex product structure is equivalent to
$\{J^{(1)},E^{(1)}_1\}$.

\

(B3') $\nabla=0$ and $\nabla'=\nabla^2$.

In this case, we have $\ggo:=\ug\bowtie\vg=\aff(\RR)\ltimes\RR^2$.
If we denote $e_i:=(e_i,0),\,f_i:=(0,f_i)$ for $i=1,2$, then the
only non zero brackets are
\[ [e_1,f_1]=\frac{1}{2}e_1,\quad [e_1,f_2]=\frac{1}{2}e_2,\quad
[e_2,f_1]=-\frac{1}{2}e_2,\quad [f_1,f_2]=f_2.\] The complex
product structure on this Lie algebra is given by $Je_i=f_i$ and
$Ee_i=e_i,\,Ef_i=-f_i$ for $i=1,2$. We will make a change of
basis, setting $v_0:=2f_1,\,v_1:=f_2,\,v_2:=2e_1,\,v_3:=-e_2$.
Thus,
\[ [v_0,v_1]=2v_1,\quad [v_0,v_2]=-v_2,\quad [v_0,v_3]=v_3,\quad [v_1,v_2]=v_3, \]
\[ Jv_0=-v_2,\quad Jv_1=v_3, \]
and then $\ggo\cong\ggo^h_2$. The eigenspaces corresponding to
$E$, in this new basis, are given by
\[ \ggo_+=\text{span}\{v_2,\,v_3\},\quad\ggo_-=\text{span}\{v_0,\,v_1\}.\]
This complex product structure is equivalent to
$\{J^{(2)},E^{(2)}_{\theta,0}=E^{(2)}_{\theta,1}\}$ with
$\theta=\pi$.

\

(B4') $\nabla=\nabla^0$ and $\nabla'=\nabla^2$.

We are looking for a linear isomorphism
$\varphi:\aff(\RR)\longrightarrow\RR^2$ compatible with $\nabla$
and $\nabla'$ and also with $\omega$ and $\omega'$. Again,
$\varphi$ must be of the form $\varphi=\begin{pmatrix} a & 0 \cr b
& d \cr \end{pmatrix}$, with $ad=1$ (in the ordered bases
$\{e_1,e_2\},\,\{f_1,f_2\}$). Hence, we have a hypersymplectic
structure on the bicrossproduct Lie algebra
$\ggo:=\ug\bowtie\vg=\RR^2\bowtie\aff(\RR)$. Let us denote
$e_i:=(e_i,0),\,f_i:=(0,f_i)$ for $i=1,2$; the only non zero
brackets are
\[ [e_1,f_1]=\frac{1}{2}e_1-abe_2+d^2f_2,\quad [e_1,f_2]=\frac{1}{2}a^2e_2,\quad
[e_2,f_1]=-\frac{1}{2}e_2,\quad [f_1,f_2]=f_2. \]
The complex product structure on this Lie algebra is given by
\[ Je_1=af_1+bf_2,\quad Je_2=df_2,\quad Jf_1=-de_1+be_2,\quad Jf_2=-ae_2 \]
and $Ee_i=e_i,\,Ef_i=-f_i$ for $i=1,2$. We will make now a change
of basis, setting $v_0:=-\frac{4}{3}d^2e_2+2f_1,\quad v_1:=f_2
,\quad v_2:=2de_1-2be_2+\frac{4}{3}d^3f_2,\quad v_3:=-ae_2$. Thus,
\[ [v_0,v_1]=2v_1,\quad [v_0,v_2]=-v_2,\quad [v_0,v_3]=v_3,\quad [v_1,v_2]=v_3, \]
\[ Jv_0=-v_2,\quad Jv_1=v_3, \]
and then $\ggo\cong\ggo^h_2$. The eigenspaces corresponding to
$E$, in this new basis, are given by
\[ \ggo_+=\text{span}\left\{-v_2+\frac{4}{3}d^3v_1,\,v_3 \right\},\quad
\ggo_-=\text{span}\left\{v_0-\frac{4}{3}d^3v_3,\,v_1\right\}.\]
This complex product structure is equivalent to
$\{J^{(2)},E^{(2)}_1\}$.

\

{\bf Case (C): $\ug=\aff(\RR)$ and $\vg=\aff(\RR)$.}

We will use $\ug=\text{span}\{e_1,e_2\}$ with $[e_1,e_2]=e_2$ and
$\vg=\text{span}\{f_1,f_2\}$ with $[f_1,f_2]=f_2$; the symplectic
forms are $\omega=e^1\wedge e^2$ and $\omega'=f^1\wedge f^2$.
Clearly, none of the connections on $\ug$ or $\vg$ may be zero,
since in that case the Lie algebra would turn out to be abelian.

\

(C1) $\nabla=\nabla^1$ and $\nabla'=\nabla^1$.

We are looking for a linear isomorphism
$\varphi:\ug\longrightarrow\vg$ compatible with $\nabla$ and
$\nabla'$ and also with $\omega$ and $\omega'$. After lengthy
computations, we obtain that $\varphi$ must be of the form
$\varphi=\begin{pmatrix} a & 0 \cr b & d \cr \end{pmatrix}$, with
$ad=1$ (in the ordered bases $\{e_1,e_2\},\,\{f_1,f_2\}$). Hence,
we have a hypersymplectic structure on the bicrossproduct Lie
algebra $\ggo:=\ug\bowtie\vg=\aff(\RR)\bowtie\aff(\RR)$. Let us
describe this Lie algebra. Let us denote
$e_i:=(e_i,0),\,f_i:=(0,f_i)$ for $i=1,2$; the only non zero
brackets are
\begin{gather*}
 [e_1,e_2]=e_2,\quad [e_1,f_1]=e_1-2abe_2-f_1-2bdf_2,\\
 [e_1,f_2]=f_2,\quad [e_2,f_1]=-e_2,\quad [f_1,f_2]=f_2.
\end{gather*}
The complex product structure on this Lie algebra is given by
\[ Je_1=af_1+bf_2,\quad Je_2=df_2,\quad Jf_1=-de_1+be_2,\quad Jf_2=-ae_2 \]
and $Ee_i=e_i,\,Ef_i=-f_i$ for $i=1,2$. We will make a change of basis, setting
\begin{gather*}
v_0:=-\frac{1}{a^2+1}\left(e_1+a^2f_1\right),\quad
v_1:=\frac{1}{a^2+1}\left(a(e_1-f_1)-a^2be_2-bf_2\right),\\
v_2:=ae_2,\quad v_3:=f_2.
\end{gather*}
Thus, we have that
\[ [v_0,v_1]=v_1,\quad [v_0,v_2]=-v_2,\quad [v_0,v_3]=-v_3,\]
\[ Jv_0=v_1,\quad Jv_2=v_3,\]
and then $\ggo\cong\ggo^h_1$. The eigenspaces corresponding to
$E$, in this new basis, are given by
\[ \ggo_+=\text{span}\left\{v_0-av_1-\frac{ab}{a^2+1}v_3,\,v_2\right\},\quad \ggo_-=
\text{span}\left\{v_0+dv_1+\frac{b}{a^2+1}v_2,\,v_3\right\}.\]
This complex product structure is equivalent to either
$\{J^{(1)},E^{(1)}_{\theta,0}\}$ if $b=0$ or
$\{J^{(1)},E^{(1)}_{\theta,1}\}$ if $b\neq 0$, in both cases with
$\cos\theta=\frac{1-a^2}{1+a^2},\, \sin\theta=-\frac{2a}{1+a^2}$.
Note that $\theta\neq 0$ and $\theta\neq\pi$ (since $a\neq 0$).

\

(C2) $\nabla=\nabla^1$ and $\nabla'=\nabla^2$ or $\nabla=\nabla^2$ and $\nabla'=\nabla^1$.

In these cases there does not exist any
$\varphi:\ug\longrightarrow\vg$ compatible with $\nabla$ and
$\nabla'$.

\

(C3) $\nabla=\nabla^2$ and $\nabla'=\nabla^2$.

We are looking for a linear isomorphism
$\varphi:\ug\longrightarrow\vg$ compatible with $\nabla$ and
$\nabla'$ and also with $\omega$ and $\omega'$. After lengthy
computations, we obtain that $\varphi$ must be of the form
$\varphi=\begin{pmatrix} a & 0 \cr b & d \cr \end{pmatrix}$, with
$ad=1$ (in the ordered bases $\{e_1,e_2\},\,\{f_1,f_2\}$). Hence,
we have a hypersymplectic structure on the bicrossproduct Lie
algebra $\ggo:=\ug\bowtie\vg=\aff(\RR)\bowtie\aff(\RR)$. Let us
describe this Lie algebra. Let us denote
$e_i:=(e_i,0),\,f_i:=(0,f_i)$ for $i=1,2$; the only non zero
brackets are
\begin{gather*}
[e_1,e_2]=e_2,\quad [e_1,f_1]=\frac{1}{2}e_1-abe_2-\frac{1}{2}f_1-bdf_2,\\
[e_1,f_2]=\frac{1}{2}(a^2e_2+f_2),\quad [e_2,f_1]=-\frac{1}{2}(e_2+d^2f_2),\quad
[f_1,f_2]=f_2.
\end{gather*}
The complex product structure on this Lie algebra is given by
\[ Je_1=af_1+bf_2,\quad Je_2=df_2,\quad Jf_1=-de_1+be_2,\quad Jf_2=-ae_2 \]
and $Ee_i=e_i,\,Ef_i=-f_i$ for $i=1,2$. We will make a change of basis, setting
\begin{gather*}
v_0:=\frac{2}{a^2+1}\left(e_1-\frac{2a^3b}{3(a^2+1)}e_2+a^2f_1+
\frac{2ab}{3(a^2+1)}f_2\right),\\
v_2:=\frac{2}{a^2+1}\left(a(e_1-f_1)-\frac{a^2b(3a^2+1)}{3(a^2+1)}e_2-
\frac{b(a^2+3)}{3(a^2+1)}f_2\right),\\
v_1:=e_2+f_2,\quad v_3:=-ae_2+df_2.
\end{gather*}
Thus, we have that
\[ [v_0,v_1]=2v_1,\quad [v_0,v_2]=-v_2,\quad [v_0,v_3]=v_3,\quad [v_1,v_2]=v_3, \]
\[ Jv_0=-v_2,\quad Jv_1=v_3, \]
and then $\ggo\cong\ggo^h_2$. The eigenspaces corresponding to
$E$, in this new basis, are given by
\begin{gather*}
\ggo_+=\text{span}\left\{v_0+av_2+\frac{2a^2b}{3(a^2+1)}v_3,\,v_1-av_3\right\},\\
\ggo_-=\text{span}\left\{v_0-dv_2-\frac{2ab}{3(a^2+1)}v_1,\,av_1+v_3\right\}.
\end{gather*}
This complex product structure is equivalent to either
$\{J^{(2)},E^{(2)}_{\theta,0}\}$ if $b=0$ or
$\{J^{(2)},E^{(2)}_{\theta,1}\}$ if $b\neq 0$, in both cases with
$\cos\theta=\frac{1-a^2}{1+a^2},\; \sin\theta=-\frac{2a}{1+a^2}$.
Note that $\theta\neq 0$ and $\theta\neq\pi$ (since $a\neq 0$).
\end{proof}

\

\begin{rem}
Complex structures on $4$-dimensional solvable Lie algebras were
classified in \cite{O} and \cite{Sn}. The Lie algebra $\ggo^h_0$
lies in the class $S1$ of \cite{Sn}, the algebra $\ggo^h_1$ lies
in the class $A2\,(\lambda=-1)$ of \cite{O} and finally,
$\ggo^h_2$ is in the class $H5\,(\lambda_1=2,\,\lambda_2=-1)$ of
\cite{O}. The first two Lie algebras carry only one complex
structure, up to equivalence, and they coincide with the complex
structures $J^{(0)}$ and $J^{(1)}$ constructed in Theorem
\ref{four}, respectively. In contrast, $\ggo^h_2$ carries two non
equivalent complex structures: one of them coincides with the
complex structure $J^{(2)}$ from Theorem \ref{four}, while the
other one cannot be part of any hypersymplectic structure on
$\ggo^h_2$.
\end{rem}

\

\section{Hypersymplectic metrics on the associated Lie groups}

We will determine now the hypersymplectic metrics on the Lie
algebras described in Theorem~\ref{four}. We will show in the
following lemma that given a complex product structure on a
4-dimensional Lie algebra, there is only one compatible metric, up
to a non zero constant. A proof of this lemma can also be found in \cite{BV}.

\begin{lem}\label{homotecia}
Let $\{J,E\}$ be a complex product structure on a 4-dimensional
Lie algebra. If $g$ and $h$ are two metrics on $\ggo$ compatible
with $\{J,E\}$, then there exists $\lambda\in\RR\setminus\{0\}$
such that $h=\lambda g$.
\end{lem}

\begin{proof}
Since $g$ and $h$ are non degenerate, there exists a linear
isomorphism $T$ of $\ggo$ such that $h(x,y)=g(Tx,y)$ for all
$x,y\in\ggo$. Let $\lambda$ be an eigenvalue of $T$, which is a
non zero real number, and let $V_\lambda$ denote the corresponding
eigenspace. Using equations (\ref{metric}), we obtain that $T$
commutes with $J$ and $E$; therefore, $V_\lambda$ is invariant by
$J$ and $E$ and then the dimension of $V_\lambda$ is even.

Let us suppose that $\dim V_\lambda=2$. As $V_\lambda$ is a proper
subspace of $\ggo$, $T$ has another eigenvalue $\mu\neq\lambda$
($\mu\in\RR\setminus\{0\}$) with $V_\mu$, the eigenspace
corresponding to $\mu$, also invariant by $J$ and $E$; hence
$\ggo=V_\lambda\oplus V_\mu$. The metric $g$ is degenerate on both
$V_\lambda$ and $V_\mu$ since otherwise the dimension of each of
these subspaces would be a multiple of 4; in fact $g=0$ on
$V_\lambda$ and on $V_\mu$.

There exist $x\in V_\lambda,\,y\in V_\mu$ such that $g(x,y)\neq
0$, since $g$ is non degenerate. We compute
\[ h(x,y)=g(Tx,y)=\lambda g(x,y) \]
and also
\[ h(x,y)=h(y,x)=g(Ty,x)=\mu g(x,y).\]
Since $g(x,y)\neq 0$, we have that $\lambda=\mu$, a contradiction.
Therefore $\dim V_\lambda=4$ and the lemma follows.
\end{proof}

\begin{rem}
\end{rem}

\

\subsection{The Lie algebra $\RR^4$}

In this case, the corresponding simply connected abelian Lie group
is the pseudo-Riemmanian manifold $\RR^{2,2}$, that is, $\RR^4$
with the neutral metric $g=\dif t^2+\dif x^2-\dif y^2-\dif z^2$,
where $t,x,y,z$ are the global canonical coordinates on $\RR^4$.
This metric is complete and flat.

\

\subsection{The Lie algebra $\ggo_0^h$}
Let $G_0^h$ denote the simply connected Lie group corresponding to
$\ggo_0^h$. It is well known that $G_0^h$ is diffeomorphic to
$\RR^4$ and, by standard computations, we can find global
coordinates $t,\,x,\,y,\,z$ on $G_0^h$ such that the
left-invariants $1$-forms $\{v^0,v^1,v^2,v^3\}$ are given by
\begin{align*}
v^0 & =\dif t, \\
v^1 & =\dif x,\\
v^2 & =\dif y,\\
v^3 & =-x\dif y+\dif z.
\end{align*}

\noi Let us consider the complex product structure
$\{J^{(0)},\,E^{(0)}_{\theta}\}$ as given in Theorem~\ref{four}.
The subalgebras $\ggo_+$ and $\ggo_-$ corresponding to
$E^{(0)}_{\theta}$ are given by
\[ \ggo_+=\text{span}\{\ctt v_3+\stt v_0,v_1\}, \quad
\ggo_-=\text{span}\{-\stt v_3+\ctt v_0,v_2\}. \] Let us denote
$U_\theta=\ctt v_3+\stt v_0$ and $V_\theta=-\stt v_3+\ctt v_0$.
Every hypersymplectic metric on $\ggo_1^h$ corresponding to this
complex product structure is homothetic to
\[{\small g_{\theta}=\begin{pmatrix} 0 & 0 & 0 & -1 \cr 0 & 0 & 1 & 0 \cr 0 & 1 & 0 & 0 \cr
-1 & 0 & 0 & 0 \end{pmatrix} }\] in the ordered basis
$\{U_{\theta},v_1,V_{\theta},v_2\}$, due to Lemma \ref{homotecia}.
If $\{v^0,v^1,v^2,v^3\}$ is the dual basis of
$\{v_0,v_1,v_2,v_3\}$, then $g_{\theta,0}$ can be written as
\begin{multline*}
g_\theta=-(\ctt v^3+\stt v^0)\otimes v^2+v^1\otimes(-\stt v^3+\ctt v^0)+\\
(-\stt v^3+\ctt v^0)\otimes v^1-v^2\otimes(\ctt v^3+\stt v^0),\end{multline*}
or, equivalently,
\[ g_\theta=-(\ctt v^3+\stt v^0)\cdot v^2+v^1\cdot(-\stt v^3+\ctt v^0)\]
where $\cdot$ denotes the symmetric product of $1$-forms. Hence,
the left-invariant metric $g_{\theta}$ on $G_0^h$ is given in
terms of the global coordinates by
\begin{multline*}
g_\theta=  -\ctt x^2\dif y^2+\ctt x\dif y\dif z-
            \ctt \dif z^2 +\stt x\dif t\dif y-\\
            \stt\dif t\dif z+\stt x \dif x\dif y-
            \stt \dif x\dif z+\ctt \dif t \dif x.
\end{multline*}
The torsion-free connection $\ncp$ on $\ggo_0^h$ associated to
$\{J^{(0)},\,E^{(0)}_{\theta}\}$ (or the Levi-Civita connection of
$g_\theta$) is easy to compute, and we can readily verify that
this connection is flat. Also, using equation (\ref{complete}), we
obtain that this connection is complete. Hence, the metrics
$g_\theta$ on $G_0^h$ are all flat and complete, and thus these
metrics are all isometric to the canonical neutral metric on
$\RR^4$, even though the hypersymplectic structures are not
equivalent.

\

\subsection{The Lie algebra $\ggo_1^h$}
Let $G_1^h$ denote the simply connected Lie group corresponding to
$\ggo_1^h$. It is well known that $G_1^h$ is diffeomorphic to
$\RR^4$ and, by standard computations, we can find global
coordinates $t,\,x,\,y,\,z$ on $G_1^h$ such that the
left-invariants $1$-forms $\{v^0,v^1,v^2,v^3\}$ are given by
\begin{align*}
v^0 & =\dif t, \\
v^1 & =e^{-t} \dif x,\\
v^2 & =e^t\dif y,\\
v^3 & =e^t\dif z.
\end{align*}

\noi $\ri$ Let us fix firstly the complex product structure
$\{J^{(1)},\,E^{(1)}_{\theta,0}\}$ as given in Theorem~\ref{four}.
The subalgebras $\ggo_+$ and $\ggo_-$ corresponding to
$E^{(1)}_{\theta,0}$ are given by
\[ \ggo_+=\text{span}\{\ctt v_0+\stt v_1,v_2\}, \quad
\ggo_-=\text{span}\{-\stt v_0+\ctt v_1,v_3\}. \] Let us denote
$U_{\theta}=\ctt v_0+\stt v_1$ and $V_{\theta}=-\stt v_0+\ctt
v_1$. Note that $JU_{\theta}=V_{\theta}$ and $Jv_2=v_3$. Every
hypersymplectic metric on $\ggo_1^h$ corresponding to this complex
product structure is homothetic to
\[ {\small g_{\theta,0}=\begin{pmatrix} 0 & 0 & 0 & -1 \cr 0 & 0 & 1 & 0 \cr
0 & 1 & 0 & 0 \cr -1 & 0 & 0 & 0 \end{pmatrix} }\] in the ordered
basis $\{U_{\theta},v_2,V_{\theta},v_3\}$, due to Lemma
\ref{homotecia}. If $\{v^0,v^1,v^2,v^3\}$ is the dual basis of
$\{v_0,v_1,v_2,v_3\}$, then $g_{\theta,0}$ can be written as
\begin{multline*}
g_{\theta,0}=-(\ctt v^0+\stt v^1)\otimes v^3+v^2\otimes(-\stt v^0+\ctt v^1)+ \\
(-\stt v^0+\ctt v^1)\otimes v^2-v^3\otimes(\ctt v^0+\stt v^1)
\end{multline*}
or, equivalently,
\[ g_{\theta,0}=-(\ctt v^0+\stt v^1)\cdot v^3 + (-\stt v^0+\ctt v^1)\cdot v^2 \]
where $\cdot$ denotes the symmetric product of $1$-forms. Hence,
the left-invariant metric $g_{\theta,0}$ on $G_1^h$ is given in
terms of the global coordinates by
\[ g_{\theta,0}=-\ctt e^t\dif t\dif z-\stt\dif x\dif z-\stt e^t\dif t\dif y+
\ctt\dif x\dif y.\]

The connection $\ncp=\nabla^{g_{\theta,0}}$ on $\ggo_1^h$ can be
explicitly computed, and we can deduce from this computation that
it is flat. However, this connection cannot be complete. Indeed,
if it were complete, then its restrictions to the eigenspaces of
$E^{(1)}_{\theta,0}$ should be complete. But at least one of these
eigenspaces is isomorphic to $\aff(\RR)$, and we know from
Proposition \ref{equiv-aff} that any flat torsion-free connection
on $\aff(\RR)$ compatible with a symplectic form is not complete,
a contradiction.

\vs

\noi $\rii$ Let us fix now the complex product structure
$\{J^{(1)},\,E^{(1)}_{\theta,1}\}$ as given in Theorem~\ref{four}.
The subalgebras $\ggo_+$ and $\ggo_-$ corresponding to
$E^{(1)}_{\theta,1}$ are given by
\begin{gather*} \ggo_+=\text{span}\{\ctt v_0+\stt v_1+\ctt v_3,v_2\},\\
\ggo_-=\text{span}\{-\stt v_0+\ctt v_1-\ctt v_2,v_3\}. \end{gather*}
Let us denote $U_{\theta}=\ctt v_0+\stt v_1+\ctt v_3$ and
$V_{\theta}=-\stt v_0+\ctt v_1-\ctt v_2$. Note that
$JU_{\theta}=V_{\theta}$ and $Jv_2=v_3$. Every hypersymplectic
metric on $\ggo_1^h$ corresponding to this complex product
structure is homothetic to
\[{\small  g_{\theta,1}=\begin{pmatrix} 0 & 0 & 0 & -1 \cr 0 & 0 & 1 & 0 \cr 0 & 1 & 0 & 0
\cr -1 & 0 & 0 & 0 \end{pmatrix} }\]
in the ordered basis $\{U_{\theta},v_2,V_{\theta},v_3\}$, because
of Lemma \ref{homotecia}. If $\{v^0,v^1,v^2,v^3\}$ is the dual
basis of $\{v_0,v_1,v_2,v_3\}$, then $g_{\theta,1}$ can be written
as
\[ g_{\theta,1}=-(\ctt v^0+\stt v^1+\ctt v^3)\cdot v^3 +
(-\stt v^0+\ctt v^1-\ctt v^2)\cdot  v^2 \]
where $\cdot$ denotes the symmetric product of $1$-forms. Hence,
the left-invariant metric $g_{\theta,1}$ on $G_1^h$ is given in
terms of the global coordinates by
\begin{multline*}
g_{\theta,1}=\ctt e^t\dif t\dif z+\stt\dif x\dif z+\ctt e^{2t}\dif z^2+ \\
\stt e^t\dif t\dif y-\ctt\dif x\dif y+\ctt e^{2t}\dif y^2.
\end{multline*}

The connection $\ncp=\nabla^{g_{\theta,1}}$ on $\ggo_1^h$ can be
explicitly computed, and we can deduce from this computation that
\[{\small  R(U_\theta,V_\theta)=\begin{pmatrix} 0&0&0&0 \cr 6\cos(\theta/2)&0&0&0\cr
0&0&0&0 \cr 0&0&6\cos(\theta/2)&0 \end{pmatrix}}\] in the ordered
basis $\{U_\theta,v_2,V_\theta,v_3\}$ and is zero for the other possibilities. Hence, $g_{\theta,1}$ is
flat if and only if $\theta=\pi$. As in the previous case, the
metrics $g_{\theta,1}$ are not complete.

\vs

\noi $\riii$ Let us consider finally the complex product structure
$\{J^{(1)},E^{(1)}_1\}$ as given in Theorem~\ref{four}. The
subalgebras $\ggo_+$ and $\ggo_-$ corresponding to $E^{(1)}_1$ are
\[ \ggo_+=\text{span}\{v_1+v_3,v_2\},\quad \ggo_-=\text{span}\{-v_0-v_2,v_3\}. \]
Every hypersymplectic metric on $\ggo_1^h$ corresponding to this
complex product structure is homothetic to
\[ g_{1}=-(v^1+v^3)\cdot v^3-(v^0+v^2)\cdot v^2,\]
which gives rise to a left-invariant metric on $G_1^h$, given by
\[ g_1=\dif x\dif z+e^{2t}\dif z^2+e^t\dif t\dif y+ e^{2t}\dif y^2.\]

It can be shown that the metric $g_1$ is flat and not complete.

\

\subsection{The Lie algebra $\ggo_2^h$}
Let $G_2^h$ denote the simply connected Lie group corresponding to
$\ggo_2^h$. It is well known that $G_2^h$ is diffeomorphic to
$\RR^4$ and, by standard computations, we can find global
coordinates $t,\,x,\,y,\,z$ on $G_2^h$ such that the
left-invariants $1$-forms $\{v^0,v^1,v^2,v^3\}$ are given by
\begin{align*}
v^0 & =\dif t, \\
v^1 & =e^{-2t} \dif x,\\
v^2 & =e^t\dif y,\\
v^3 & =e^{-t}(\dif z-\ft x\dif y+\ft y\dif x).
\end{align*}

\noi $\ri$ Let us fix firstly the complex product structure
$\{J^{(2)},\,E^{(2)}_{\theta,0}\}$ as given in Theorem~\ref{four}.
The subalgebras $\ggo_+$ and $\ggo_-$ corresponding to
$E^{(2)}_{\theta,0}$ are given by
\begin{gather*} \ggo_+=\text{span}\{\ctt v_0+\stt v_2,\ctt v_1-\stt v_3\}, \\
\ggo_-=\text{span}\{\stt v_0-\ctt v_2,\stt v_1+\ctt v_3\}.
\end{gather*} Let us denote $U_{\theta}=\ctt v_0+\stt
v_1,\,{\tilde U}_\theta=\ctt v_1-\stt v_3$ and $V_{\theta}=-\stt
v_0+\ctt v_1,\, {\tilde V}_\theta=\stt v_1+\ctt v_3$. Note that
$JU_{\theta}=V_{\theta}$ and $J{\tilde U}_\theta={\tilde
V}_\theta$. Every hypersymplectic metric on $\ggo_2^h$
corresponding to this complex product structure is homothetic to
\[{\small  g_{\theta,0}=\begin{pmatrix} 0 & 0 & 0 & -1 \cr 0 & 0 & 1 & 0 \cr
0 & 1 & 0 & 0 \cr -1 & 0 & 0 & 0 \end{pmatrix} }\] in the ordered
basis $\{U_{\theta},{\tilde U}_\theta,V_{\theta},{\tilde
V}_\theta\}$, due to Lemma \ref{homotecia}. If
$\{v^0,v^1,v^2,v^3\}$ is the dual basis of $\{v_0,v_1,v_2,v_3\}$,
then $g_{\theta,0}$ can be written as
\begin{multline*} 
g_{\theta,0}=-(\ctt v^0+\stt v^2)\cdot(\stt v^1+\ctt v^3)+\\
 (\ctt v^1-\stt v^3)\cdot(\stt v^1+\ctt v^3)  \end{multline*}
where $\cdot$ denotes the symmetric product of $1$-forms. Hence,
the left-invariant metric $g_{\theta,0}$ on $G_2^h$ is given in
terms of the global coordinates by
\[ g_{\theta,0}=e^{-t}\dif t(\dif z-\ft x\dif y+\ft y\dif x)+e^{-t}\dif x\dif y. \]
Note that $g_{\theta,0}$ does not depend on $\theta$. It can be
shown that this metric is flat and not complete.

\vs

\noi $\rii$ Let us fix now the complex product structure
$\{J^{(2)},\,E^{(2)}_{\theta,2}\}$ as given in Theorem~\ref{four}.
The subalgebras $\ggo_+$ and $\ggo_-$ corresponding to
$E^{(2)}_{\theta,1}$ are given by
\begin{gather*} \ggo_+=\text{span}\{\ctt v_0+\stt v_2+\ctt v_3,\ctt v_1-\stt v_3\}\\
\ggo_-=\text{span}\{\stt v_0-\ctt v_2-\ctt v_1,\stt v_1+\ctt
v_3\}. \end{gather*} Let us denote $U_{\theta}=\ctt v_0+\stt
v_2+\ctt v_3,\,{\tilde U}_\theta=\ctt v_1-\stt v_3$ and
$V_{\theta}=\stt v_0-\ctt v_2-\ctt v_1,\, {\tilde V}_\theta=\stt
v_1+\ctt v_3$. Note that $JU_{\theta}=V_{\theta}$ and $J{\tilde
U}_\theta={\tilde V}_\theta$. Every hypersymplectic metric on
$\ggo_2^h$ corresponding to this complex product structure is
homothetic to
\[ {\small g_{\theta,1}=\begin{pmatrix} 0 & 0 & 0 & -1 \cr 0 & 0 & 1 & 0 \cr 0 & 1 & 0 & 0
\cr -1 & 0 & 0 & 0 \end{pmatrix} }\] in the ordered basis
$\{U_{\theta},{\tilde U}_{\theta},V_{\theta},{\tilde
V}_{\theta}\}$, because of Lemma \ref{homotecia}. If
$\{v^0,v^1,v^2,v^3\}$ is the dual basis of $\{v_0,v_1,v_2,v_3\}$,
then $g_{\theta,1}$ can be written as
\begin{multline*}
g_{\theta,1}=-(\ctt v^0+\stt v^2+\ctt v^3)\cdot(\stt v^1+\ctt v^3)+ \\
 (\stt v^0-\ctt v^2-\ctt v^1)\cdot(\ctt v^1-\stt v^3)  \end{multline*}
where $\cdot$ denotes the symmetric product of $1$-forms. Hence,
the left-invariant metric $g_{\theta,1}$ on $G_2^h$ is given in
terms of the global coordinates by
\begin{multline*}
g_{\theta,1}=e^{-t}\dif t(\dif z-\ft x\dif y+\ft y\dif x)+e^{-t}\dif x\dif y+
\cos^2_{\theta/2}e^{-4t}\dif x^2+\\
\cos^2_{\theta/2}e^{-2t}(\dif z-\ft x\dif y+\ft y\dif x)^2.
\end{multline*}

It can be shown that the curvature of the connection
$\ncp=\nabla^{g_{\theta,1}}$ is given by
\[ {\small R(U_\theta,V_\theta)=\begin{pmatrix} 0&0&0&0 \cr 6\cos^2(\theta/2)&0&0&0\cr
0&0&0&0 \cr 0&0&6\cos^2(\theta/2)&0 \end{pmatrix}}\] in the
ordered basis $\{U_\theta,{\tilde U}_\theta,V_\theta,{\tilde
V}_\theta\}$ and is zero for the other possibilities. Hence,
$g_{\theta,1}$ is flat if and only if $\theta=\pi$. As in previous
cases, the metrics $g_{\theta,1}$ are not complete.

\vs

\noi $\riii$ Let us consider finally the complex product structure
$\{J^{(2)},E^{(2)}_1\}$ as given in Theorem~\ref{four}. The
subalgebras $\ggo_+$ and $\ggo_-$ corresponding to $E^{(2)}_1$ are
\[ \ggo_+=\text{span}\{v_1+v_2,v_3\},\quad \ggo_-=\text{span}\{v_0+v_3,-v_1\}. \]
Every hypersymplectic metric on $\ggo_2^h$ corresponding to this
complex product structure is homothetic to
\[ g_{1}=(v^1+v^2)\cdot v^1+(v^0+v^3)\cdot v^3,\]
which gives rise to a left-invariant metric on $G_2^h$, given by
\[ g_1=e^{-t}\dif x\dif y+e^{-4t}\dif x^2 +e^{-t}\dif t(\dif z-\ft x\dif y+
\ft y\dif x)+  e^{-2t}(\dif z-\ft x\dif y+\ft y\dif x)^2.\]

This metric is flat and not complete.

\

\end{document}